# SAMPLED-DATA STABILIZATION OF NONLINEAR DELAY SYSTEMS WITH A COMPACT ABSORBING SET

**Iasson Karafyllis[*] and Miroslav Krstic[**]**

[*]Dept. of Mathematics, National Technical University of Athens, Zografou Campus, 15780, Athens, Greece, email: iasonkar@central.ntua.gr

[**]Dept. of Mechanical and Aerospace Eng., University of California, San Diego, La Jolla, CA 92093-0411, U.S.A., email: krstic@ucsd.edu

**Abstract**
We present a methodology for the global sampled-data stabilization of systems with a compact absorbing set and input/measurement delays. The methodology is based on the Inter-Sample-Predictor, Observer, Predictor, Delay-Free Controller (ISP-O-P-DFC) scheme and the stabilization is robust to perturbations of the sampling schedule. The obtained results are novel even for the delay-free case.

**Keywords:** sampled-data control, nonlinear systems, delay systems.

## 1. Introduction

Achieving stabilization by sampled-data output feedback and ensuring robustness to perturbations in the sampling schedule are central challenges in nonlinear control over networks, where the simultaneous presence of asynchrony and (measurement or input) delays creates important problems (see [7,8,9,24,25,27,28,29,31]). Almost all available results rely on delay-dependent conditions for the existence of stabilizing feedback and in most cases the stability domain depends on the sampling interval/ delay. Predictive feedback seems to be the only possible choice for handling large delays (see [3,4,5,11,12,13,14,15,18,19,26]). Global stabilization of control systems with large delays by means of sampled-data output feedback with positive sampling rate remains a challenging problem. There are few results on the global stabilization of systems with input applied with Zero-Order-Hold and sampled measurements which do not coincide with the state vector (output measurement) even in the delay-free case; see [12,2,26,6,20]. The existing results either exploit the linear structure or a global Lipschitz property. In general, global results for sampled-data output feedback control of delayed systems are limited; see also [22] for results with sufficiently small delays.

The present work provides global stabilization results for a class of nonlinear systems: systems with a compact absorbing set. More specifically, we consider nonlinear systems of the form

$$\dot{x}(t) = f(x(t), u(t-\tau)), x \in \Re^n, u \in U \quad (1.1)$$

where $U \subseteq \Re^m$ is a non-empty compact set with $0 \in U$, $\tau \geq 0$ is the input delay and $f: \Re^n \times \Re^m \to \Re^n$ is a smooth vector field with $f(0,0) = 0$. The measurements are sampled and the output is given by

$$y(\tau_i) = h(x(\tau_i - r)) \quad (1.2)$$

where $h: \Re^n \to \Re^k$ is a smooth mapping with $h(0) = 0$ and $r \geq 0$ is the measurement delay. The class of nonlinear systems of the form (1.1), (1.2), with a compact absorbing set has been studied in [1,10,13]. Here, we extend the ideas in [13] to the case where the input is applied with Zero-Order-Hold (ZOH) and we use the Inter-Sample-Predictor, Observer, Predictor, Delay-Free Controller (ISP-O-P-DFC) control scheme. The ISP-O-P-DFC control scheme has long been in use for linear systems [21,23,30,32]. The main idea of the control scheme is the use of an inter-sample predictor of the (not available) continuous output signal. The observer uses the approximation of the continuous output signal and provides an estimate of the delayed state vector, which is subsequently fed to an approximate predictor: the predictor provides an estimation of the future value of the state vector. Finally, the estimation of the future value of the state is used by the delay-free controller and the control action is applied with Zero-Order-Hold. A major difference with [13] (except of the fact that [13] considered continuously applied input) is the predictor: the approximate predictor used in the present work is the repeated explicit Euler numerical scheme for the control system (1.1). The prediction scheme was used in [14,15] and can be easily implemented in computer software (since in the present work the applied time step is constant).

Our main result (Theorem 2.2) provides explicit formulas for global stabilizers, which are robust with respect to perturbations of the sampling schedule. Moreover, Theorem 2.2 can be also applied to the case where the sampling times do not necessarily coincide with the times that input changes. This feature is important for network systems and is rare in the sampled-data control literature (usually the sample-and-hold case is studied). The state is driven to the equilibrium at an exponential rate. The result of Theorem 2.2 is novel even for the delay-free case $r = \tau = 0$. Corollary 2.3 presents a specialization of the result to the delay-free case. See also [17] for semi-global results in the delay-free case based on sampled-data dynamic output feedback.

The structure of the present work is as follows: Section 2 is devoted to the presentation of the basic assumptions for nonlinear systems with a compact absorbing set and the statement of the main results. The proof of the main result is provided in Section 3, where additional lemmas are stated and utilized. An illustrative example is shown in Section 4, where the proposed control scheme is applied. The concluding remarks are provided in Section 5. Finally, the Appendix contains the proofs of all auxiliary lemmas used in Section 3.

**Notation.** Throughout this paper, we adopt the following notations:

* $\Re_+ := [0,+\infty)$. A partition of $\Re_+$ is an increasing sequence $\{\tau_i\}_{i=0}^{\infty}$ with $\tau_0 = 0$ and $\lim_{i \to \infty} \tau_i = +\infty$.

* Let $I \subseteq \Re_+ := [0,+\infty)$ be an interval. By $L^{\infty}(I;U)$, we denote the space of measurable and essentially bounded functions $u(\cdot)$ defined on $I$ and taking values in $U \subseteq \Re^m$. Let $A \subseteq \Re^n$ be an open set. By $C^0(A;\Omega)$, we denote the class of continuous functions on $A$, which take values in $\Omega \subseteq \Re^m$. By $C^k(A;\Omega)$, where $k \geq 1$ is an integer, we denote the class of functions on $A \subseteq \Re^n$ with continuous derivatives of order $k$, which take values in $\Omega \subseteq \Re^m$. For a function $V = (V_1,...,V_k)' \in C^1(A;\Re^k)$, the gradient of $V$ at $x \in A \subseteq \Re^n$, denoted by $\nabla V(x)$, is a matrix with $k$ rows; its $i$-th row is the row vector $\left[ \frac{\partial V_i}{\partial x_1}(x) \quad ... \quad \frac{\partial V_i}{\partial x_n}(x) \right]$ for $i = 1,...,k$.

* For a vector $x \in \Re^n$, we denote by $x'$ its transpose and by $|x|$ its Euclidean norm. $A' \in \Re^{n \times m}$ denotes the transpose of the matrix $A \in \Re^{m \times n}$ and $|A|$ denotes the induced norm of the matrix $A \in \Re^{m \times n}$, i.e., $|A| = \sup\{|Ax| : x \in \Re^m, |x| = 1\}$. $I \in \Re^{n \times n}$ denotes the unit matrix.

* A function $V : \Re^n \to \Re_+$ is called positive definite if $V(0) = 0$ and $V(x) > 0$ for all $x \neq 0$. A function $V : \Re^n \to \Re_+$ is called radially unbounded if the sets $\{x \in \Re^n : V(x) \leq M\}$ are either empty or bounded for all $M \geq 0$.

* The class of functions $K_{\infty}$ is the class of strictly increasing, continuous functions $a : \Re_+ \to \Re_+$ with $a(0) = 0$ and $\lim_{s \to +\infty} a(s) = +\infty$. For $x \in \Re$, $[x]$ denotes the integer part of $x \in \Re$.

* For $u : [a-r,b) \to U$, where $U \subseteq \Re^m$, $b > a$ and $r > 0$, $u_t : [-r,0] \to U$ for $t \in [a,b)$ denotes the $r$-"history" of $u$, i.e., the function defined by $(u_t)(\theta) = u(t+\theta)$ for $\theta \in [-r,0]$ and $\tilde{u}_t : [-r,0) \to U$ for $t \in [a,b]$ denotes the $r$-"open history" of $u$, i.e., the function defined by $(u_t)(\theta) = u(t+\theta)$ for $\theta \in [-r,0)$. For a bounded function $u : [-r,0] \to U$ (or $u : [-r,0) \to U$), $\|u\|$ denotes the norm $\|u\| = \sup_{-r \leq \theta \leq 0}(|u(\theta)|)$ (or $\|u\| = \sup_{-r \leq \theta < 0}(|u(\theta)|)$).

## 2. Problem Description and Main Result

Our first assumption for system (1.1) guarantees that there exists a compact set which is robustly globally asymptotically stable. We call the compact set "absorbing" because the solution "is absorbed" in the set after an initial transient period.

**(H1)** *There exist a radially unbounded (but not necessarily positive definite) function $V \in C^2(\Re^n; \Re_+)$, a positive definite function $W \in C^1(\Re^n; \Re_+)$ and a constant $R > 0$ such that the following inequality holds for all $(x,u) \in \Re^n \times U$ with $V(x) \geq R$*

$$\nabla V(x) f(x,u) \leq -W(x) \quad (2.1)$$

*Moreover, the set $S_1 = \{x \in \Re^n : V(x) \leq R\}$ contains a neighborhood of $0 \in \Re^n$.*

Indeed, assumption (H1) guarantees that for every initial condition $x(0) \in \Re^n$ and for every measurable and essentially bounded input $u : \Re_+ \to U$ the solution $x(t)$ of (1.1) enters the compact set $S_1 = \{x \in \Re^n : V(x) \leq R\}$ after a finite transient period, i.e., there exists $T \in C^0(\Re^n; \Re_+)$ such that $x(t) \in S$, for all $t \geq T(x(0))$. Moreover, notice that the compact set $S_1 = \{x \in \Re^n : V(x) \leq R\}$ is positively invariant. This fact is guaranteed by the following lemma which is an extension of Theorem 5.1 in [16] (page 211). The proof of the following lemma can be found in [10].

**Lemma 2.1:** *Consider system (1.1) under assumption (H1). There exists $T \in C^0(\Re^n; \Re_+)$ such that for every $x_0 \in \Re^n$ and for every measurable and essentially bounded input $u : [-\tau, +\infty) \to U$ the solution $x(t) \in \Re^n$ of (1.1) with initial condition $x(0) = x_0$ and corresponding to input $u : [-\tau, +\infty) \to U$ satisfies $V(x(t)) \leq \max(V(x_0), R)$ for all $t \geq 0$ and $V(x(t)) \leq R$ for all $t \geq T(x_0)$.*

Our second assumption guarantees that we are in a position to construct an appropriate local exponential stabilizer for the delay-free version system (1.1), i.e., system (1.1) with $\tau = 0$.

**(H2)** *There exist a positive definite function $P \in C^2(\Re^n; \Re_+)$, constants $\mu, K_1 > 0$ with $K_1 |x|^2 \leq P(x)$ for all $x \in \Re^n$ with $V(x) \leq R$ and a locally Lipschitz mapping $k : \Re^n \to U$ with $k(0) = 0$ such that the following inequality holds*

$$\nabla P(x) f(x, k(x)) \leq -2\mu |x|^2, \text{ for all } x \in \Re^n \text{ with } V(x) \leq R \quad (2.2)$$

Our third assumption guarantees that we are in a position to construct an appropriate local exponential observer for the delay-free system (1.1), (1.2) with $r = \tau = 0$.

**(H3)** *There exist a symmetric and positive definite matrix $Q \in \Re^{n \times n}$, constants $\omega > 0$, $b > R$ and a matrix $L \in \Re^{n \times k}$ such that the following inequality holds*

$$(z-x)' Q \big( f(z,u) + L(h(z) - h(x)) - f(x,u) \big) \leq -\omega |z-x|^2,$$
$$\text{for all } u \in U, z, x \in \Re^n \text{ with } V(z) \leq b \text{ and } V(x) \leq R \quad (2.3)$$

Our final assumption is a technical assumption that enables us to construct a dynamic feedback stabilizer for system (1.1), (1.2). Similar assumptions have been used in [1,10,13].

**(H4)** *There exist constants $c \in (0,1)$, $R \leq a < b$ such that the following inequality holds:*

$$\nabla V(z)(f(z,u) + L(h(z) - h(x))) \leq -W(z) + (1-c) |\nabla V(z)|^2 \frac{(z-x)' Q \big( f(z,u) + L(h(z) - h(x)) - f(x,u) \big)}{\nabla V(z) Q(z-x)}$$
$$\text{for all } u \in U, z, x \in \Re^n \text{ with } a < V(z) \leq b, \nabla V(z) Q(z-x) < 0 \text{ and } V(x) \leq R \quad (2.4)$$

Assumption (H4) implies restrictions on the dynamics of the local observer, which was introduced by assumption by (H3). Notice that the left hand side of inequality (2.4) is the time derivative of the function $V(z(t))$ along the trajectories of the local observer $\dot{z} = f(z,u) + L(h(z) - h(x))$ with input $x \in \Re^n$. Therefore, assumption (H4) imposes an upper bound on the time derivative of the function $V(z(t))$ along the trajectories of the local observer $\dot{z} = f(z,u) + L(h(z) - h(x))$ with input $x \in \Re^n$ for certain regions of the state space: the solution of the local observer is not allowed to "grow too fast".

Assumption (H4) is needed for a specific reason. Lemma 2.1 implies that the states of any successful observer for system (1.1), (1.2) with $r = \tau = 0$ must be driven in a compact set after a transient period. Therefore, the design of an observer with a compact absorbing set is desired. Assumption (H4) is a sufficient condition that allows us to design a global observer with a compact absorbing set (expressed by the sublevel sets of $V$) which coincides with the local observer $\dot{z} = f(z,u) + L(h(z) - h(x))$ on an appropriate neighborhood of the equilibrium. In order to achieve the design of such an observer, we need to impose bounds on the "growth" of the trajectories of the local observer $\dot{z} = f(z,u) + L(h(z) - h(x))$ with input $x \in \Re^n$ for certain regions of the state space. However, it should be noticed that (2.4) does not exclude the possibility of having a positive time derivative of the function $V(z(t))$ along the trajectories of the local observer $\dot{z} = f(z,u) + L(h(z) - h(x))$.

We are now in a position to state the main result of the present work.

**Theorem 2.2:** *Consider system (1.1), (1.2) under assumptions (H1-4). Define:*

$$\hat{k}(z,y,u) := L(h(z) - y), \text{ for all } (z,y,u) \in \Re^n \times \Re^k \times U \text{ with } V(z) \le R \quad (2.5)$$

$$\hat{k}(z,y,u) := L(h(z) - y) - \frac{\varphi(z,y,u)}{|\nabla V(z)|^2}(\nabla V(z))', \text{ for all } (z,y,u) \in \Re^n \times \Re^k \times U \text{ with } V(z) > R \quad (2.6)$$

*where $\varphi : \Re^n \times \Re^k \times \Re^m \to \Re_+$ is defined by*

$$\varphi(z,y,u) := \max\big(0, \nabla V(z) f(z,u) + W(z) + p(V(z)) \nabla V(z) L(h(z) - y)\big) \quad (2.7)$$

*and $p : \Re_+ \to [0,1]$ is an arbitrary locally Lipschitz function that satisfies $p(s) = 1$ for all $s \ge b$ and $p(s) = 0$ for all $s \le a$. Let $N > 0$ be an integer and define the mapping:*

$$\Phi_N : \Re^n \times L^\infty\big([-r-\tau,0);U\big) \to \Re^n \quad (2.8)$$

*which maps $(x_0, u) \in \Re^n \times L^\infty\big([-r-\tau,0);U\big)$ to the vector $\Phi_N(x_0,u) := x_N \in \Re^n$, where $x_i \in \Re^n$ ($i=1,...,N$) are vectors given by the recursive formula:*

$$x_{i+1} = x_i + \int_{ih}^{(i+1)h} f(x_i, u(s-r-\tau)) ds, \text{ for } i = 0,...,N-1 \quad (2.9)$$

*where $h := (\tau + r)/N$. Then for sufficiently small constants $T_s > 0$, $T_H > 0$ and for sufficiently large integer $N > 0$ there exist a locally Lipschitz function $C \in K_\infty$ and a constant $\sigma > 0$ such that for every partition $\{\tau_i\}_{i=0}^\infty$ of $\Re_+$ with $\sup_{i \ge 0}(\tau_{i+1} - \tau_i) \le T_s$, $z_0 \in \Re^n$, $x_0 \in C^0\big([-r,0];\Re^n\big)$, $u_0 \in L^\infty\big([-r-\tau,0);U\big)$, the solution of (1.1), (1.2) with*

$$\dot{z}(t) = f(z(t), u(t-r-\tau)) + \hat{k}(z(t), w(t), u(t-r-\tau)), \text{ for } t \ge 0 \text{ a.e.} \quad (2.10)$$

$$\dot{w}(t) = \nabla h(z(t)) f(z(t), u(t-r-\tau)), \text{ for } t \in [\tau_i, \tau_{i+1}) \text{ a.e. and for all integers } i \ge 0 \quad (2.11)$$

$$w(\tau_i) = y(\tau_i), \text{ for all integers } i \ge 0 \quad (2.12)$$

$$u(t) = k\big(\Phi_N(z(jT_H), \breve{u}_{jT_H})\big), \text{ for all } t \in [jT_H, (j+1)T_H) \text{ and for all integers } j \ge 0 \quad (2.13)$$

*initial condition $z(0) = z_0$, $x(\theta) = x_0(\theta)$ for $\theta \in [-r,0]$, $u(\theta) = u_0(\theta)$ for $\theta \in [-r-\tau,0)$, exists and satisfies the following estimate for all $t \ge 0$:*

$$\|x_t\| + |z(t)| + \|\breve{u}_t\| \le \exp(-\sigma t) C\big(\|x_0\| + |z_0| + \|\breve{u}_0\|\big) \quad (2.14)$$

The result of Theorem 2.2 is novel even for the delay-free case $r = \tau = 0$. Indeed, one can repeat the proof of Theorem 2.2 and obtain the following corollary (its proof is omitted due to the similarity with the proof of Theorem 2.2).

**Corollary 2.3:** *Consider system (1.1), (1.2) under assumptions (H1-4) with $r = \tau = 0$. Let $\varphi : \Re^n \times \Re^k \times \Re^m \to \Re_+$ be defined by (2.7) for an arbitrary locally Lipschitz function $p : \Re_+ \to [0,1]$ that satisfies $p(s) = 1$ for all $s \geq b$ and $p(s) = 0$ for all $s \leq a$ and let $\hat{k}(z, y, u)$ be defined for all $(z, y, u) \in \Re^n \times \Re^k \times U$ by (2.5), (2.6). Then for sufficiently small constants $T_s > 0$, $T_H > 0$ there exist a locally Lipschitz function $C \in K_\infty$ and a constant $\sigma > 0$ such that for every partition $\{\tau_i\}_{i=0}^\infty$ of $\Re_+$ with $\sup_{i \geq 0}(\tau_{i+1} - \tau_i) \leq T_s$, $z_0 \in \Re^n$, $x_0 \in \Re^n$, the solution of (1.1), (1.2) with*

$$\dot{z}(t) = f(z(t), u(t)) + \hat{k}(z(t), w(t), u(t)), \text{ for } t \geq 0 \text{ a.e.} \tag{2.15}$$

$$\dot{w}(t) = \nabla h(z(t)) f(z(t), u(t)), \text{ for } t \in [\tau_i, \tau_{i+1}) \text{ a.e. and for all integers } i \geq 0 \tag{2.16}$$

$$w(\tau_i) = y(\tau_i), \text{ for all integers } i \geq 0 \tag{2.17}$$

$$u(t) = k(z(jT_H)), \text{ for all } t \in [jT_H, (j+1)T_H) \text{ and for all integers } j \geq 0 \tag{2.18}$$

*initial condition $z(0) = z_0$, $x(0) = x_0$, exists and satisfies the following estimate for all $t \geq 0$:*

$$|x(t)| + |z(t)| \leq \exp(-\sigma t) C(|x_0| + |z_0|) \tag{2.19}$$

**Remark 2.4:** (a) It should be emphasized that estimate (2.14) guarantees (i) local exponential stabilization (since $C \in K_\infty$ is locally Lipschitz) and (ii) global $K$-exponential stabilization.
(b) The approximate predictor mapping given (2.8), (2.9) is the repeated explicit Euler numerical scheme for the control system (1.1). It can be easily implemented in computer software.
(c) The proof of Theorem 2.2 is constructive. Therefore, estimates of the size of the constants $T_s > 0$, $T_H > 0$, $N > 0$ and $\sigma > 0$ are provided. However, the estimates are conservative.
(d) The construction of the controller (2.10), (2.11), (2.12), (2.13) (or (2.15), (2.16), (2.17), (2.18)) is based on the local controller provided by assumption (H2) and the local observer provided by assumption (H3). However, the results of Theorem 2.2 and Corollary 2.3 are global.
(e) The construction of the controller (2.10), (2.11), (2.12), (2.13) is based on the Inter-Sample-Predictor, Observer, Predictor, Delay-Free Controller (ISP-O-P-DFC) control scheme. Namely, (2.11), (2.12) is the inter-sample predictor of the continuous output signal and is fed to the observer (2.10). The observer estimate is fed to the approximate predictor mapping given by (2.8), (2.9), which provides the estimation $\Phi_N(z(jT_H), \breve{u}_{jT_H})$ of the future value of the state vector. Finally, the estimation $\Phi_N(z(jT_H), \breve{u}_{jT_H})$ of the future value of the state is fed to the controller (2.13) and the control action is applied with ZOH. It is important to notice that if measurement delay is present (i.e., if $r > 0$), then the observer provides an estimate of the delayed state vector $x(t - r)$.

## 3. Proof of Main Result

Define the sets:

$$S_1 := \{x \in \Re^n : V(x) \leq R\}, \quad S_2 := \{x \in \Re^n : V(x) \leq b\} \tag{3.1}$$

Notice that since $V \in C^2(\Re^n; \Re_+)$ is radially unbounded, it follows that the sets defined in (3.1) are compact sets.

The proof of the main result requires a number of technical lemmas. The first technical lemma provides an estimate for the observation error.

**Lemma 3.1:** *Let $\sigma > 0$, $T_s > 0$ be sufficiently small constants and let $\{\tau_i\}_{i=0}^\infty$ be a partition of $\Re_+$ with $\sup_{i \geq 0}(\tau_{i+1} - \tau_i) \leq T_s$. Then there exists a constant $M_1 > 0$ such that every solution of (1.1), (1.2), (2.10), (2.11), (2.12), corresponding to (arbitrary) input $u \in L^\infty([-r-\tau, +\infty); U)$ and satisfying $x(t-r) \in S_1$, $z(t) \in S_2$ for all $t \geq \tau_l$, where $l > 0$ is an integer with $\tau_l \geq r$, also satisfies the following inequality for all $t \geq \tau_l$:*

$$\sup_{\tau_l \leq s \leq t} (\exp(\sigma s) |z(s) - x(s - r)|) \leq M_1 \exp(\sigma \tau_l) |z(\tau_l) - x(\tau_l - r)| \tag{3.2}$$

The second technical lemma provides an estimate for the state.

**Lemma 3.2:** *Let $\sigma > 0$, $T_H > 0$ be sufficiently small constants. Then there exist constants $M_2, M_3 > 0$ such that every solution of (1.1), (1.2), corresponding to (arbitrary) input $u \in L^\infty([-\tau, +\infty); U)$ and satisfying $x(t) \in S_1$ for all $t \geq \tau + jT_H$, where $j \geq 0$ is an integer, also satisfies the following inequality for all $t \geq jT_H$:*

$$\sup_{jT_H \leq s \leq t} \left( |x(s+\tau)| \exp(\sigma s) \right) \leq M_2 \exp(\sigma j T_H) |x(jT_H + \tau)| + M_3 \sup_{jT_H \leq s \leq t} \left( \left| u(s) - k\left( x\left( \tau + \left[\frac{s}{T_H}\right] T_H \right) \right) \right| \exp(\sigma s) \right) \quad (3.3)$$

The third technical lemma provides an estimate for the prediction error.

**Lemma 3.3:** *There exists an integer $N^* > 0$ and a constant $M_4 > 0$ such that for every $N \geq N^*$ for every $x_0 \in S_2$ and for every measurable and essentially bounded input $u \in L^\infty([-r-\tau, 0); U)$ the following estimates hold for the solution $x(t)$ of (1.1) with initial condition $x(-r) = x_0$, corresponding to (arbitrary) input $u \in L^\infty([-r-\tau, 0); U)$:*

$$|x(\tau) - \Phi_N(x_0, u)| \leq \frac{M_4}{N} (|x_0| + \|u\|) \quad (3.4)$$

$$x_i \in S_2, \text{ for all } i = 0, 1, \ldots, N \quad (3.5)$$

*where $\|u\| = \sup_{-r-\tau \leq s < 0} (|u(s)|)$ and $x_i \in \Re^n$ ($i = 1, \ldots, N$) are vectors given by the recursive formula (2.9).*

The fourth technical lemma the three previous lemmas and provides an estimate for the norm of the solution of the closed-loop system (1.1), (1.2), (2.10), (2.11), (2.12), (2.13).

**Lemma 3.4:** *Let $\sigma > 0$, $T_s > 0$, $T_H > 0$ be sufficiently small constants and let $N \geq 1$ be a sufficiently large integer. Let $\{\tau_i\}_{i=0}^\infty$ be a partition of $\Re_+$ with $\sup_{i \geq 0}(\tau_{i+1} - \tau_i) \leq T_s$. Then there exists a constant $G > 0$ such that every solution of (1.1), (1.2), (2.10), (2.11), (2.12), (2.13) satisfying $x(t-r-T_H) \in S_1$, $z(t) \in S_2$ for all $t \geq jT_H$, where $j \geq 0$ is an integer, also satisfies the following inequality for all $t \geq 0$:*

$$\left( \|\breve{u}_t\| + \|x_t\| + |z(t)| \right) \exp(\sigma t) \leq G \left( \sup_{-r-\tau \leq s \leq jT_H} (|u(s)|) + \sup_{0 \leq s \leq r+\tau+jT_H+T_s} (|z(s)|) + \sup_{-r \leq s \leq jT_H+T_s+\tau+r} (|x(s)|) \right) \quad (3.6)$$

We are now ready to provide the proof of Theorem 2.2.

**Proof of Theorem 2.2:** We first notice that the following inequality holds for all $(z, w, u) \in \Re^n \times \Re^k \times U$ with $V(z) \geq b$:

$$\nabla V(z)(f(z,u) + \hat{k}(z,w,u)) \leq -W(z) \quad (3.7)$$

Definition (2.6) implies $\nabla V(z)(f(z,u) + \hat{k}(z,w,u)) = \nabla V(z)(f(z,u) + L(h(z)-w)) - \varphi(z,w,u)$. By distinguishing the cases $\nabla V(z)f(z,u) + W(z) + \nabla V(z)L(h(z)-w) \leq 0$ and $\nabla V(z)f(z,u) + W(z) + \nabla V(z)L(h(z)-w) > 0$, using definition (2.7) and noticing that $p(V(z)) = 1$ we conclude that (3.7) holds.

Let $\sigma > 0$, $T_s > 0$, $T_H > 0$ be sufficiently small constants and let $N \geq 1$ be a sufficiently large integer so that Lemma 3.4 holds. Let $\{\tau_i\}_{i=0}^\infty$ be a partition of $\Re_+$ with $\sup_{i \geq 0}(\tau_{i+1} - \tau_i) \leq T_s$ and let $x_0 \in C^0([-r, 0]; \Re^n)$, $z_0 \in \Re^n$, $u_0 \in L^\infty([-r-\tau, 0); U)$ be given. We will show first that the solution of (1.1), (1.2), (2.10), (2.11), (2.12), (2.13), with initial condition $z(0) = z_0$, $x(\theta) = x_0(\theta)$ for $\theta \in [-r, 0]$, $u(\theta) = u_0(\theta)$ for $\theta \in [-r-\tau, 0)$ exists for all $t \geq 0$ and is unique.

We first make the following claim.

**Claim 1:** Suppose that $x(t)$ is defined on $[-r, \tau_{i+1}]$, $u(t)$ is defined on $[-r-\tau, \tau_{i+1})$ and that $z(t)$ is defined on $[0, \tau_i]$. Then $z(t)$ is defined on $[0, \tau_{i+1}]$.

Standard results in ordinary differential equations guarantee that the system
$$\dot{z}(t) = f(z(t), u(t-r-\tau)) + \hat{k}(z(t), w(t), u(t-r-\tau))$$
$$\dot{w}(t) = \nabla h(z(t)) f(z(t), u(t-r-\tau))$$
(3.8)

has a local solution defined on $[\tau_i, \tilde{t})$ for some $\tilde{t} \in (\tau_i, \tau_{i+1}]$. By virtue of (3.7) and Lemma 2.1, it follows that the solution of (3.8) satisfies the following estimate:
$$V(z(t)) \leq \max(V(z_0), b) \tag{3.9}$$
for all $t \geq 0$ for which the solution of (3.8) exists. Define the non-decreasing function:
$$\Omega(s) := \max\{|\nabla h(z) f(z, u)| : (z, u) \in \Re^n \times U, V(z) \leq s\}, \text{ for all } s \geq \min(V(z) : z \in \Re^n) \tag{3.10}$$
which is well-defined by virtue of the facts that $U \subseteq \Re^m$ is compact and $V \in C^2(\Re^n; \Re_+)$ is a radially unbounded function. It follows from definition (3.10) and inequality (3.9), that the solution of (3.8) satisfies the following estimate for all $t \in [\tau_i, \tilde{t})$:
$$|w(t)| \leq |w(\tau_i)| + T_s \Omega(\max(V(z_0), b)) \tag{3.11}$$
A standard contradiction argument shows that $z(t)$ is defined on $[0, \tau_{i+1}]$.

The second claim guarantees existence/uniqueness of solutions for all $t \geq 0$. It is an application of the method of steps.

**Claim 2:** $\breve{u}_t$, $x_t$, $z(t)$ are uniquely determined for all $t \in [0, jT_H]$, where $j \in Z_+$.

The claim is proved by induction. First we notice that the claim holds for $j = 0$. Next, we show that if the claim holds for some $j \in Z_+$ then the claim holds for $j+1$. Indeed, (2.13) guarantees that $\breve{u}_t$ is uniquely determined for all $t \in (jT_H, (j+1)T_H]$. It follows from Lemma 2.1 that $x_t$ is uniquely determined for all $t \in (jT_H, (j+1)T_H]$. Since the set $(jT_H, (j+1)T_H] \cap \{\tau_i\}_{i=0}^\infty$ is either empty or finite, Claim 1 implies that we are in a position determine uniquely $z(t)$ for all $t \in (jT_H, (j+1)T_H]$. Thus $\breve{u}_t$, $x_t$, $z(t)$ are uniquely determined for all $t \in [0, (j+1)T_H]$, where $j \in Z_+$.

Lemma 2.1 in conjunction with (2.1) and (3.7) implies there exists $T \in C^0(\Re^n; \Re_+)$ such that the inequalities $V(x(t)) \leq \max(V(x_0), R)$ and (3.9) hold for all $t \geq 0$ and
$$V(x(t)) \leq R \text{ for all } t \geq T(x_0(0)) \text{ and } V(z(t)) \leq b \text{ for all } t \geq T(z_0) \tag{3.12}$$
Indeed, the above conclusions for $V(x(t))$ are direct consequences of Lemma 2.1. The above conclusions for $V(z(t))$ are consequences of Lemma 2.1 applied to system (2.10) with $(w, u)$ as inputs. Inequalities (3.12) and definitions (3.1) show that
$$x(t) \in S_1, \ z(t) \in S_2, \text{ for all } t \geq \max(T(x_0(0)), T(z_0)) \tag{3.13}$$
Let $j \geq 0$ be the smallest integer so that $jT_H \geq r + T_H + \max(T(x_0(0)), T(z_0))$. Then (3.13) in conjunction with Lemma 3.4 implies the existence of a constant $G \geq 0$ such that (3.6) holds.

Since $f : \Re^n \times \Re^m \to \Re^n$, $k : \Re^n \to U$, $h : \Re^n \to \Re^k$, $\nabla h : \Re^n \to \Re^{k \times n}$, are locally Lipschitz mappings with $f(0,0) = 0$, $k(0) = 0$, $h(0) = 0$ and since $U \subset \Re^m$ is compact, there exists a continuous, non-decreasing function $L : \Re_+ \to [1, +\infty)$ such that:
$$\begin{aligned}|f(x,u)| + |\nabla h(x) f(x,u)| &\leq L(|x|)(|x| + |u|) \\ |h(x)| + |k(x)| &\leq L(|x|)|x|\end{aligned}, \text{ for all } x \in \Re^n, u \in U \tag{3.14}$$
Moreover, taking into account definitions (2.5), (2.6), (2.7) and inequalities (3.14), we are in a position to conclude that there exists a continuous, non-decreasing function $\hat{L} : \Re_+ \to [1, +\infty)$ such that:
$$|f(z,u) + \hat{k}(z,w,u)| \leq \hat{L}(|z|)(|z| + |u| + |w|), \text{ for all } z \in \Re^n, w \in \Re^k, u \in U \tag{3.15}$$

Furthermore, using induction, (3.14), definitions (2.8), (2.9) and the fact that $U \subset \Re^m$ is compact, we are in a position to guarantee that there exists a continuous, non-decreasing function $\tilde{L}: \Re_+ \to [1,+\infty)$ such that:

$$|\Phi_N(x,u)| \leq \tilde{L}(|x|)(|x|+\|u\|), \text{ for all } (x,u) \in \Re^n \times L^\infty([-r-\tau,0];U) \quad (3.16)$$

Using inequalities (3.14), (3.15), (3.16), we show that the following claim holds.

**Claim 3:** There exists a sequence of non-decreasing functions $\Gamma_i: \Re_+ \to \Re_+$ with $\Gamma_i(s) \leq \Gamma_{i+1}(s)$ for all $s \geq 0$ and for all integers $i \geq 0$ such that the following inequality holds for every integer $i \geq 0$:

$$\sup_{0 \leq t \leq iT_H} (\|\breve{u}_t\| + \|x_t\| + |z(t)|) \leq (\|u_0\| + \|x_0\| + |z_0|)\Gamma_i(\|x_0\| + |z_0|) \quad (3.17)$$

We construct the sequence inductively. Inequality (3.17) holds for $i = 0$ with the function $\Gamma_0(s) \equiv 1$. In order to show the claim, we assume that there exists an integer $i \geq 0$ and a non-decreasing function $\Gamma_i: \Re_+ \to \Re_+$ such that (3.17) holds. We next show that there exists a non-decreasing function $\Gamma_{i+1}: \Re_+ \to \Re_+$ with $\Gamma_i(s) \leq \Gamma_{i+1}(s)$ for all $s \geq 0$ such that (3.17) holds with $i \geq 0$ replaced by $i+1$.

Using (2.13), (3.16) and (3.17), we get for $t \in [iT_H, (i+1)T_H)$:

$$|u(t)| \leq |\Phi_N(z(iT_H), \breve{u}_{iT_H})| \leq \tilde{L}(|z(iT_H)|)(|z(iT_H)| + \|\breve{u}_{iT_H}\|)$$

$$\leq \tilde{L}((\|u_0\| + \|x_0\| + |z_0|)\Gamma_i(\|x_0\| + |z_0|))\Gamma_i(\|x_0\| + |z_0|)(\|u_0\| + \|x_0\| + |z_0|)$$

Using (3.17), the above inequality and the fact that $U \subset \Re^m$ is compact, we obtain the existence of a non-decreasing function $Z_1: \Re_+ \to \Re_+$ such that:

$$\sup_{0 \leq t \leq (i+1)T_H} (\|\breve{u}_t\|) \leq (\|u_0\| + \|x_0\| + |z_0|)Z_1(\|x_0\| + |z_0|) \quad (3.18)$$

Next, define the following family of sets for all $p \geq 0$:

$$S(p) := \{x \in \Re^n : V(x) \leq b + \max\{V(\xi) : \xi \in \Re^n, |\xi| \leq p\}\} \quad (3.19)$$

Notice that by virtue of assumption (H1) the above sets are compact for each $p \geq 0$ and that $S(p_1) \subseteq S(p_2)$ for every $p_1, p_2 \geq 0$ with $p_1 \leq p_2$. Define the non-decreasing function for all $p \geq 0$:

$$\varphi(p) := \max_{x \in S(p)} (|x|) \quad (3.20)$$

Applying the inequality $|x(t)| \leq |x(iT_H)| + \int_{iT_H}^{t} |f(x(s), u(s-\tau))| ds$ for the solution $x(t)$ of (1.1) with $t \in [iT_H, (i+1)T_H]$ and using (3.14), (3.19), (3.20) in conjunction with Lemma 2.1 and the Gronwall-Bellman Lemma, we obtain:

$$|x(t)| \leq \exp(L(\varphi(|x(0)|))T_H)\left(|x(iT_H)| + \sup_{0 \leq s \leq (i+1)T_H} (\|\breve{u}_s\|)\right), \text{ for all } t \in [iT_H, (i+1)T_H]$$

Using (3.17), (3.18) and the above inequality, we obtain the existence of a non-decreasing function $Z_2: \Re_+ \to \Re_+$ such that:

$$\sup_{0 \leq t \leq (i+1)T_H} (\|\breve{u}_t\| + \|x_t\|) \leq (\|u_0\| + \|x_0\| + |z_0|)Z_2(\|x_0\| + |z_0|) \quad (3.21)$$

Let $t_i \in \{\tau_j\}_{j=0}^{\infty}$ be the largest sampling time with $t_i \leq iT_H$. Using (2.11), (2.12), (3.14), (3.9), definitions (3.19), (3.20), we get for all $t \in [t_i, (i+1)T_H]$:

$$|w(t)| \leq A + L(\varphi(|z(0)|))\int_{t_i}^{t} |z(s)| ds + L(\varphi(|z(0)|))\int_{t_i}^{t} |u(s)| ds \quad (3.22)$$

where $A := L\left(\sup_{0 \leq t \leq (i+1)T_H}(\|x_t\|)\right) \sup_{0 \leq t \leq (i+1)T_H}(\|x_t\|)$. Using (2.10), (3.15), (3.9), definitions (3.19), (3.20), we get for all $t \in [t_i, (i+1)T_H]$:

$$|z(t)| \leq |z(t_i)| + \hat{L}(\varphi(|z(0)|))\int_{t_i}^{t}|z(s)|ds + \hat{L}(\varphi(|z(0)|))\int_{t_i}^{t}|u(s)|ds + \hat{L}(\varphi(|z(0)|))\int_{t_i}^{t}|w(s)|ds$$

$$\leq |z(t_i)| + \hat{L}(\varphi(|z(0)|))(t-t_i)A + \hat{L}(\varphi(|z(0)|))(1 + L(\varphi(|z(0)|))(t-t_i))\int_{t_i}^{t}|z(s)|ds \quad (3.23)$$

$$+ \hat{L}(\varphi(|z(0)|))(1 + L(\varphi(|z(0)|))(t-t_i))\int_{t_i}^{t}|u(s)|ds$$

Using the fact that $t_i \in \{\tau_j\}_{j=0}^{\infty}$ is the largest sampling time with $t_i \leq iT_H$, in conjunction with $\sup_{i \geq 0}(\tau_{i+1} - \tau_i) \leq T_s$, we obtain that $t_i \geq iT_H - T_s$. Therefore, we obtain from (3.23) for all $t \in [t_i, (i+1)T_H]$:

$$|z(t)| \leq B + \tilde{\varphi}(|z(0)|)\int_{t_i}^{t}|z(s)|ds \quad (3.24)$$

where

$$\tilde{\varphi}(s) := \hat{L}(\varphi(s))(1 + L(\varphi(s))(T_H + T_s))$$
$$B := |z(t_i)| + \hat{L}(\varphi(|z(0)|))(T_H + T_s)A + (T_H + T_s)\tilde{\varphi}(|z(0)|)\sup_{0 \leq s \leq (i+1)T_H}(\|\tilde{u}_s\|) \quad (3.25)$$

Using the Gronwall-Bellman Lemma in conjunction with (3.24) and the fact that $t_i \geq iT_H - T_s$, we get for all $t \in [t_i, (i+1)T_H]$:

$$|z(t)| \leq \exp(\tilde{\varphi}(|z(0)|)(T_H + T_s))B \quad (3.26)$$

Using (3.17), (3.21), (3.25), (3.26), the fact that $A := L\left(\sup_{0 \leq t \leq (i+1)T_H}(\|x_t\|)\right)\sup_{0 \leq t \leq (i+1)T_H}(\|x_t\|)$, we are in a position to conclude that there exists a non-decreasing function $\Gamma_{i+1} : \Re_+ \to \Re_+$ such that (3.17) holds with $i \geq 0$ replaced by $i+1$.

Since $T \in C^0(\Re^n; \Re_+)$ is continuous, there exists a constant $\Omega$ and a function $\kappa \in K_\infty$ such that $T(x) \leq \Omega + \kappa(|x|)$ for all $x \in \Re^n$. Since $j \geq 0$ is the smallest integer so that $jT_H \geq r + T_H + \max(T(x_0(0)), T(z_0))$, it follows that $iT_H \geq r + jT_H + \tau + T_s$ for $i = \psi(\|x_0\| + |z_0|) = 3 + \left[\dfrac{2r + \tau + T_s + 2\kappa(\|x_0\| + |z_0|)}{T_H}\right]$. Combining (3.6) with (3.17), we obtain the following estimate for all $t \geq 0$:

$$\|x_t\| + |z(t)| + \|\tilde{u}_t\| \leq \exp(-\sigma t)(\|x_0\| + |z_0| + \|\tilde{u}_0\|)\tilde{C}(\|x_0\| + |z_0|) \quad (3.27)$$

where $\tilde{C}(s) = G\Gamma_{\psi(s)}(s)$ for all $s \geq 0$. Since $\tilde{C} : \Re_+ \to \Re_+$ is a non-decreasing function, there exists a $C^1$ function $\hat{C} : \Re_+ \to \Re_+$ such that $\hat{C}(s) \geq \tilde{C}(s)$ for all $s \geq 0$. Inequality (2.14) is a direct consequence of (3.27) and the definition $C(s) := s\hat{C}(s)$ for all $s \geq 0$. The proof is complete. ◁

## 4. Illustrative Example

Consider the following planar nonlinear system:

$$\dot{x}_1 = \zeta x_1 - 10x_1^3 + x_2 \quad, \quad \dot{x}_2 = -\frac{13}{4}x_2 + u \quad ; \quad x = (x_1, x_2)' \in \Re^2, u \in \Re \quad (4.1)$$

where $\zeta > 0$ is a constant that satisfies the inequality

$$25001\zeta^2 + 2\zeta \leq 4 \quad (4.2)$$

with output

$$y = h(x) = x_1 \quad (4.3)$$

We show next that system (4.1) satisfies assumptions (H1), (H2), (H3) and (H4) with

$$U = \left[-50\zeta\sqrt{2}, 50\zeta\sqrt{2}\right] \quad (4.4)$$

Therefore, Theorem 2.2 can be applied to (4.1) and the system can be stabilized for arbitrary input and measurement delays by bounded feedback applied with ZOH. The reader (who is used in continuous feedback stabilization for delay-free nonlinear systems) may be surprised by the existence of an upper bound for the constant $\zeta > 0$ (see (4.2)), i.e., the linear part of system (4.1) is only weakly destabilizing. Two things must be noted at this point:

(a) The upper bound for $\zeta > 0$ in (4.2) is restrictive and can be improved considerably. However, we have given this restrictive bound for simplicity (the algebraic manipulations become easier).

(b) We intend to design a feedback law for system (4.1) that: (i) is bounded, (ii) is applied with ZOH (even though the system is not linear or globally Lipschitz), (iii) uses sampled and delayed measurements with uncertain sampling schedule, (iv) guarantees stabilization in the presence of (arbitrary) input and measurement delays, (v) guarantees local exponential stabilization and a global exponential convergence rate, and (vi) can handle sampling times which do not necessarily coincide with the times that the input changes value (i.e., it is not necessarily a sample-and-hold feedback). However, in order to achieve all the above features we have to assume that there is an upper bound for $\zeta > 0$. More specifically, the requirement of the existence of a compact absorbing set (i.e., assumption (H1)) implies that the input $u$ of system (4.1) takes values in a compact set $U$. The simultaneous requirement of having a local exponential stabilizer for (4.1) (i.e., assumption (H2)) leads to the fact that the size of $U$, the size of the set where the local exponential stabilizer works and the constant $\zeta > 0$ are related and an upper bound for $\zeta > 0$ is needed.

Assumption (H1) holds with $V(x) = (x_1^2 + x_2^2)/2$, $R = 1$ and $W(x) = V(x)/4$. Indeed, we get:

$$\nabla V(x) f(x, u) = \zeta x_1^2 - 10 x_1^4 + x_1 x_2 - \frac{13}{4} x_2^2 + x_2 u$$

Using the inequalities $x_1 x_2 \leq (x_1^2 + x_2^2)/2$, $u x_2 \leq (u^2 + x_2^2)/2$ and (4.4), we obtain for all $(x, u) \in \Re^2 \times U$:

$$\nabla V(x) f(x, u) \leq -W(x) + \left( \zeta + \frac{5}{8} - 10 x_1^2 \right) x_1^2 - \frac{17}{8} x_2^2 + 2500 \zeta^2 \quad (4.5)$$

If $10 x_1^2 > \zeta + 1$, then (4.5) implies $\nabla V(x) f(x, u) \leq -W(x) - \frac{3}{4} V(x) + 2500 \zeta^2$. By virtue of (4.2), the previous inequality directly implies (2.1) for the case $10 x_1^2 > \zeta + 1$. If $10 x_1^2 \leq \zeta + 1$, then (4.5) implies $\nabla V(x) f(x, u) \leq -W(x) + \frac{1}{10}(\zeta + 1)^2 - \frac{3}{4} V(x) + 2500 \zeta^2$. By virtue of (4.2), the previous inequality directly implies (2.1) for the case $10 x_1^2 \leq \zeta + 1$. Therefore, we conclude that (2.1) holds in every case.

Assumption (H2) holds with $P(x) = \frac{1}{2} x_1^2 + \frac{2}{\zeta(13 - 4\zeta)}(x_2 + 2\zeta x_1)^2$, $\tilde{k}(x) = -\frac{3\zeta}{4}(13 - 4\zeta) x_1 + 20 \zeta x_1^3$, $k(x) = \min\left( 50 \zeta \sqrt{2}, \max\left( -50 \zeta \sqrt{2}, \tilde{k}(x) \right) \right)$ and appropriate constants $\mu, K_1 > 0$. Notice that the fact that $P(x) = \frac{1}{2} x_1^2 + \frac{2}{\zeta(13 - 4\zeta)}(x_2 + 2\zeta x_1)^2$ is a quadratic positive definite function implies the existence of a constant $K_1 > 0$ with $K_1 |x|^2 \leq P(x)$ for all $x \in \Re^2$. Moreover, by virtue of (4.2), we get for all $x \in \Re^2$ with $V(x) = (x_1^2 + x_2^2)/2 \leq 1 = R$:

$$\left| \tilde{k}(x) \right| \leq \frac{3\zeta}{4}(13 - 4\zeta)\sqrt{2} + 40 \zeta \sqrt{2} \leq 50 \zeta \sqrt{2} \quad (4.6)$$

Therefore, the equality $k(x) = \tilde{k}(x) = -\frac{3\zeta}{4}(13 - 4\zeta) x_1 + 20 \zeta x_1^3$ holds for all $x \in \Re^2$ with $V(x) = (x_1^2 + x_2^2)/2 \leq 1 = R$. Notice that, by virtue of (4.2), the following inequality holds for all $x \in \Re^2$ with $V(x) = (x_1^2 + x_2^2)/2 \leq 1 = R$:

$$\nabla P(x) f(x, k(x)) = \nabla P(x) f(x, \tilde{k}(x)) = -\zeta x_1^2 - 10 x_1^4 - \frac{13 - 8\zeta}{a(13 - 4\zeta)}(x_2 + 2\zeta x_1)^2 \leq -2\zeta P(x) \quad (4.7)$$

Inequality (4.7) in conjunction with the fact that $P(x) = \frac{1}{2} x_1^2 + \frac{2}{\zeta(13 - 4\zeta)}(x_2 + 2\zeta x_1)^2$ is a quadratic positive definite function implies the existence of a constant $\mu > 0$ such that (2.2) holds.

Next we show that assumption (H3) holds with $Q = I \in \Re^{2\times 2}$, $L = -(2\zeta, 1)'$, $\omega = \zeta > 0$ and arbitrary constant $b > 1 = R$. Indeed, we have by virtue of (4.2), for all $(x, z, u) \in \Re^2 \times \Re^2 \times U$

$$(z-x)'Q(f(z,u) + L(h(z) - h(x))) - f(x,u)) = -\zeta(z_1 - x_1)^2 - 10(z_1^2 + z_1 x_1 + x_1^2)(z_1 - x_1)^2 - \frac{13}{4}(z_2 - x_2)^2 \leq -\zeta|z - x|^2$$

which holds because $z_1^2 + z_1 x_1 + x_1^2 \geq 0$ for all $(x_1, z_1) \in \Re^2$ and because $\zeta \leq 13/4$.

Finally, we show that assumption (H4) holds. More specifically, we show that the more demanding inequality

$$\nabla V(z)(f(z,u) + L(h(z) - h(x))) = \zeta z_1^2 - 10 z_1^4 + (z_1 + u)z_2 - (2\zeta z_1 + z_2)(z_1 - x_1) - \frac{13}{4}z_2^2 \leq -\frac{1}{8}(z_1^2 + z_2^2) = -W(z) \quad (4.8)$$

holds for all $u \in U$, $z, x \in \Re^2$ with $a < V(z)$ and $V(x) \leq 1 = R$ for sufficiently large $a \geq 1$. Therefore, (2.4) holds with arbitrary constants $c \in (0,1)$ and $a < b$. Inequality (4.8) is equivalent to the inequality

$$u z_2 + 2\zeta z_1 x_1 + x_1 z_2 \leq \left(\zeta + 10 z_1^2 - \frac{1}{8}\right)z_1^2 + \frac{25}{8} z_2^2$$ which, by virtue of (4.4) and the fact that $V(x) \leq 1 = R$, is directly implied by the inequality

$$2500 \zeta^2 + 2\zeta \sqrt{2} |z_1| + |z_2|\sqrt{2} \leq \left(\zeta + 10 z_1^2 - \frac{1}{8}\right)z_1^2 + \frac{21}{8} z_2^2 \quad (4.9)$$

Similarly using the inequalities $|z_2|\sqrt{2} \leq 1 + z_2^2/2$ and $2\zeta\sqrt{2}|z_1| \leq z_1^2 + 2\zeta^2$, we conclude that (4.9) holds provided that the following inequality holds:

$$2502 \zeta^2 \leq \left(\zeta + 10 z_1^2 - \frac{9}{8}\right)z_1^2 + \frac{17}{8} z_2^2 \quad (4.10)$$

If $10 z_1^2 > 26/8$, then (4.10) holds for $20016 \zeta^2/17 \leq z_1^2 + z_2^2$. On the other hand, if $10 z_1^2 \leq 26/8$, then (4.10) is implied by the inequality $2502\zeta^2 + \frac{17}{8} z_1^2 \leq z_1^2 + z_2^2$, which follows from the inequality $2502 \zeta^2 + \frac{221}{320} \leq z_1^2 + z_2^2$. We conclude that (4.10) holds for all $u \in U$, $z, x \in \Re^2$ provided that $V(z) = \frac{1}{2} z_1^2 + \frac{1}{2} z_2^2 \geq a = \max\left(1, \frac{10008}{17}\zeta^2, 1251\zeta^2 + \frac{221}{640}\right)$.

Let $b > a$ be an arbitrary constant and let $p : \Re_+ \to [0,1]$ be an arbitrary locally Lipschitz function that satisfies $p(s) = 1$ for all $s \geq b$ and $p(s) = 0$ for all $s \leq a$. Define

$$\hat{k}(z, y, u) := -\begin{bmatrix} 2\zeta \\ 1 \end{bmatrix}(z_1 - y), \text{ for all } (z, y, u) \in \Re^2 \times \Re \times U \text{ with } z_1^2 + z_2^2 \leq 2 \quad (4.11)$$

$$\hat{k}(z, y, u) := -\begin{bmatrix} 2\zeta \\ 1 \end{bmatrix}(z_1 - y) - \frac{\varphi(z, y, u)}{z_1^2 + z_2^2}\begin{bmatrix} z_1 \\ z_2 \end{bmatrix}, \text{ for all } (z, y, u) \in \Re^2 \times \Re \times U \text{ with } z_1^2 + z_2^2 > 2 \quad (4.12)$$

where $\varphi : \Re^2 \times \Re \times \Re \to \Re_+$ is defined by

$$\varphi(z, y, u) := \max\left(0, \left(\zeta + \frac{1}{8} - 10 z_1^2\right)z_1^2 + (z_1 + u)z_2 - \frac{25}{8} z_2^2 - p\left(\frac{z_1^2 + z_2^2}{2}\right)(2\zeta z_1 + z_2)(z_1 - y)\right) \quad (4.13)$$

Let $\tau, r \geq 0$ be arbitrary constants. Theorem 2.2 guarantees that for sufficiently small constants $T_s > 0$, $T_H > 0$ and for sufficiently large integer $N > 0$ there exist a locally Lipschitz function $C \in K_\infty$ and a constant $\sigma > 0$ such that for every partition $\{\tau_i\}_{i=0}^\infty$ of $\Re_+$ with $\sup_{i \geq 0}(\tau_{i+1} - \tau_i) \leq T_s$, $z_0 \in \Re^2$, $x_0 \in C^0([-r, 0]; \Re^2)$, $u_0 \in L^\infty([-r - \tau, 0); U)$, the solution of

$$\dot{x}_1(t) = \zeta x_1(t) - 10 x_1^3(t) + x_2(t) \; ; \; \dot{x}_2(t) = -\frac{13}{4} x_2(t) + u(t - \tau) \quad (4.14)$$

$$\dot{z}(t) = \begin{bmatrix} \zeta z_1(t) - 10 z_1^3(t) + z_2(t) \\ -\frac{13}{4} z_2(t) + u(t - r - \tau) \end{bmatrix} + \hat{k}(z(t), w(t), u(t - r - \tau)), \text{ for } t \geq 0 \text{ a.e.} \quad (4.15)$$

$$\dot{w}(t) = \zeta z_1(t) - 10 z_1^3(t) + z_2(t), \text{ for } t \in [\tau_i, \tau_{i+1}) \text{ a.e. and for all integers } i \geq 0 \quad (4.16)$$

$$w(\tau_i) = x_1(\tau_i - r), \text{ for all integers } i \geq 0 \tag{4.17}$$

$$u(t) = \min\left(50\zeta\sqrt{2}, \max\left(-50\zeta\sqrt{2}, -\frac{3\zeta}{4}(13-4\zeta)q_N^{(1)}(jT_H) + 20\zeta\left(q_N^{(1)}(jT_H)\right)^3\right)\right),$$

$$\text{for all } t \in [jT_H, (j+1)T_H) \text{ and for all integers } j \geq 0 \tag{4.18}$$

$$\begin{bmatrix} q_0^{(1)}(jT_H) \\ q_0^{(2)}(jT_H) \end{bmatrix} = z(jT_H), \begin{bmatrix} q_{i+1}^{(1)}(jT_H) \\ q_{i+1}^{(2)}(jT_H) \end{bmatrix} = \begin{bmatrix} (1+\zeta h)q_i^{(1)}(jT_H) - 10h\left(q_i^{(1)}(jT_H)\right)^3 + hq_i^{(2)}(jT_H) \\ \left(1-\frac{13h}{4}\right)q_i^{(2)}(jT_H) + \int_{ih}^{(i+1)h} u(jT_H + s - r - \tau)ds \end{bmatrix}, \text{ for } i = 0,\ldots,N-1$$

$$\tag{4.19}$$

where $h := (r+\tau)/N$, initial condition $z(0) = z_0$, $x(\theta) = x_0(\theta)$ for $\theta \in [-r,0]$, $u(\theta) = u_0(\theta)$ for $\theta \in [-r-\tau, 0)$, exists and satisfies estimate (2.14) for all $t \geq 0$.

This example shows that even if the delay-free system can be globally stabilized by a static output feedback, still an observer must be used when delays are present. The reason that forces the use of the observer is the prediction: in order to make an accurate prediction for the future value of the output, accurate estimates of the state vector are needed. Indeed, system (4.1) can be globally stabilized by the static output feedback $k(x) = \min\left(50\zeta\sqrt{2}, \max\left(-50\zeta\sqrt{2}, \tilde{k}(x)\right)\right)$ with $\tilde{k}(x) = -\frac{3\zeta}{4}(13-4\zeta)x_1 + 20\zeta x_1^3$. However, the dynamic feedback stabilizer given by (4.15), (4.16), (4.17), (4.18) and (4.19) uses the hybrid sampled-data observer (4.15), (4.16), (4.17): the observer state is used in the prediction scheme given by (4.19).

## 5. Concluding Remarks

The present work provides a methodology for the global sampled-data stabilization of systems with a compact absorbing set and input/measurement delays. The methodology is based on the ISP-O-P-DFC scheme and the stabilization is robust to perturbations of the sampling schedule. The obtained results are novel even for the delay-free case.

More remains to be done. The results can be extended (under appropriate assumptions) to the case where the absorbing set is not necessarily compact: the absorbing set can be a set where a Lipschitz inequality holds. This will be the topic of future research.

**References**

[1] Ahmed-Ali, T., I. Karafyllis and F. Lamnabhi-Lagarrigue, "Global Exponential Sampled-Data Observers for Nonlinear Systems with Delayed Measurements", *Systems and Control Letters*, 62(7), 2013, 539-549.

[2] Ahmed-Ali, T., I. Karafyllis, M. Krstic and F. Lamnabhi-Lagarrigue, "Robust Stabilization of Nonlinear Globally Lipschitz Delay Systems", to appear in the edited book entitled *Recent Results on Nonlinear Time Delayed Systems* (Eds. M. Malisoff, F. Mazenc, P. Pepe and I. Karafyllis) in the series Advances in Delays and Dynamics (ADD@S), Springer.

[3] Bekiaris-Liberis, N. and M. Krstic, "Compensation of State-Dependent Input Delay for Nonlinear Systems", *IEEE Transactions on Automatic Control*, 58, 2013, 275-289.

[4] Bekiaris-Liberis, N. and M. Krstic, *Nonlinear Control Under Nonconstant Delays*, SIAM, 2013.

[5] Bekiaris-Liberis, N. and M. Krstic, "Robustness of Nonlinear Predictor Feedback Laws to Time- and State-Dependent Delay Perturbations", *Automatica*, 49(6), 2013, 1576-1590.

[6] Castillo-Toledo, B., S. Di Gennaro and G. Sandoval Castro, "Stability Analysis for a Class of Sampled Nonlinear Systems With Time-Delay", *Proceedings of the 49th Conference on Decision and Control*, Atlanta, GA, USA, 2010, 1575-1580.

[7] Fridman, E., A. Seuret and J.-P. Richard, "Robust Sampled-Data Stabilization of Linear Systems: An Input Delay Approach", *Automatica*, 40(8), 2004, 1441-1446.

[8] Gao, H., T. Chen and J. Lam, "A New System Approach to Network-Based Control", *Automatica*, 44(1), 2008, 39-52.
[9] Heemels, M., A.R. Teel, N. van de Wouw and D. Nešić, "Networked Control Systems with Communication Constraints: Tradeoffs between Transmission Intervals, Delays and Performance", *IEEE Transactions on Automatic Control*, 55(8), 2010, 1781 - 1796.
[10] Karafyllis, I. and C. Kravaris, "Global Exponential Observers for Two Classes of Nonlinear Systems", *Systems and Control Letters*, 61(7), 2012, 797-806.
[11] Karafyllis, I. and M. Krstic, "Nonlinear Stabilization under Sampled and Delayed Measurements, and with Inputs Subject to Delay and Zero-Order Hold", *IEEE Transactions on Automatic Control*, 57(5), 2012, 1141-1154.
[12] Karafyllis, I. and M. Krstic, "Stabilization of Nonlinear Delay Systems Using Approximate Predictors and High-Gain Observers", *Automatica*, 49(12), 2013, 3623–3631.
[13] Karafyllis, I., M. Krstic, T. Ahmed-Ali and F. Lamnabhi-Lagarrigue, "Global Stabilization of Nonlinear Delay Systems with a Compact Absorbing Set", *International Journal of Control*, 87(5), 2014, 1010-1027.
[14] Karafyllis, I. and M. Krstic, "Numerical Schemes for Nonlinear Predictor Feedback", to appear in *Mathematics of Control, Signals and Systems*, 26(4), 2014, 519-546.
[15] Karafyllis, I., M. Malisoff, M. de Queiroz, M. Krstic and R. Yang, "Predictor-Based Tracking for Neuromuscular Electrical Stimulation", to appear in the *International Journal of Robust and Nonlinear Control* (see also arXiv:1310.1857 [math.OC]).
[16] Khalil, H. K., *Nonlinear Systems*, 2nd Edition, Prentice-Hall, 1996.
[17] Khalil, H. K., "Performance Recovery Under Output Feedback Sampled-Data Stabilization of a Class of Nonlinear Systems", *IEEE Transactions on Automatic Control*, 49, 2004, 2173-2184.
[18] Krstic, M., *Delay Compensation for Nonlinear, Adaptive, and PDE Systems*, Birkhäuser Boston, 2009.
[19] Krstic, M., "Input Delay Compensation for Forward Complete and Strict-Feedforward Nonlinear Systems", *IEEE Transactions on Automatic Control*, 55(2), 2010, 287-303.
[20] Lozano, R., P. Castillo, P. Garcia and A. Dzul, "Robust Prediction-Based Control for Unstable Delay Systems: Application to the Yaw Control of a Mini-Helicopter", *Automatica*, 40(4), 2004, 603-612.
[21] Lozano, R., A. Sanchez, S. Salazar-Cruz and I. Fantoni, "Discrete-Time Stabilization of Integrators in Cascade: Real-Time Stabilization of a Mini-Rotorcraft" *International Journal of Control*, 81(6), 2008, 894–904.
[22] Mazenc, F., M. Malisoff, T. Dinh, "Robustness of Nonlinear Systems With respect to Delay and Sampling of the Controls", *Automatica*, 49(6), 2013, 1925-1931.
[23] Medvedev, A. and H. Toivonen, "Continuous-Time Deadbeat Observation Problem With Application to Predictive Control of Systems With Delay", *Kybernetika*, 30(6), 1994, 669–688.
[24] Nešić, D. and A. Teel, "A Framework for Stabilization of Nonlinear Sampled-Data Systems Based on their Approximate Discrete-Time Models", *IEEE Transactions on Automatic Control*, 49(7), 2004, 1103-1122.
[25] Nešić, D., A. R. Teel and D. Carnevale, "Explicit Computation of the Sampling Period in Emulation of Controllers for Nonlinear Sampled-Data Systems", *IEEE Transactions on Automatic Control*, 54(3), 2009, 619-624.
[26] Sun, X.-M., K.-Z. Liu, C. Wen and W. Wang, "Predictive Control of Nonlinear Continuous Networked Control Systems with Large Time-Varying Transmission Delays and Transmission Protocols", *personal communication*.
[27] Tabbara, M., D. Nešić and A. R. Teel, "Networked Control Systems: Emulation Based Design", in *Networked Control Systems* (Eds. D. Liu and F.-Y. Wang) Series in Intelligent Control and Intelligent Automation, World Scientific, 2007.
[28] Tabuada, P., "Event-Triggered Real-Time Scheduling of Stabilizing Control Tasks", *IEEE Transactions on Automatic Control*, 52(9), 2007, 1680-1685.
[29] Walsh, G. C., O. Beldiman, and L. G. Bushnell, "Asymptotic Behavior of Nonlinear Networked Control Systems", *IEEE Transactions on Automatic Control*, 46(7), 2001, 1093-1097.

[30] Watanabe, K., and M. Sato, "A Predictor Control for Multivariable Systems With General Delays in Inputs and Outputs Subject to Unmeasurable Disturbances", *International Journal of Control*, 40(3), 1984, 435–448.

[31] Yu, M., L. Wang, T. Chu and F. Hao, "Stabilization of Networked Control Systems with Data Packet Dropout and Transmissions Delays: Continuous-Time Case", *European Journal of Control*, 11(1), 2005, 41-49.

[32] Zhong, Q.-C., *Robust Control of Time-Delay Systems*, Springer-Verlag, London, 2010.

# Appendix

**Proof of Lemma 3.1:** First we establish the following inequality:

$$(z-x)'Q\big(f(z,u)+\hat{k}(z,h(x),u)-f(x,u)\big) \leq -c\omega|z-x|^2, \text{ for all } (x,z,u)\in S_1\times S_2\times U \quad \text{(A.1)}$$

Notice that inequalities (2.1), (2.3) and definitions (2.5), (2.6), (2.7) imply that (A.1) holds for the case $V(z)\leq a$. Therefore, we focus on the case $a<V(z)\leq b$. Definition (2.6) gives:

$$(z-x)'Q\big(f(z,u)+\hat{k}(z,h(x),u)-f(x,u)\big) \leq (z-x)'Q\big(f(z,u)+L(h(z)-h(x))-f(x,u)\big) - \frac{\varphi(z,h(x),u)}{|\nabla V(z)|^2}\nabla V(z)Q(z-x) \quad \text{(A.2)}$$

Inequalities (2.3), (A.2) and the fact that $\varphi(z,h(x),u)\geq 0$ implies that (A.1) holds if $\nabla V(z)Q(z-x)\geq 0$. Moreover, inequalities (2.3), (A.2) show that (A.1) holds if $\varphi(z,h(x),u)=0$. It remains to consider the case $\nabla V(z)Q(z-x)<0$ and $\varphi(z,h(x),u)>0$. In this case, definition (2.7) implies $\varphi(z,h(x),u)=\nabla V(z)f(z,u)+W(z)+p(V(z))\nabla V(z)L(h(z)-h(x))>0$. Then, inequality (2.4) gives:

$$\varphi(z,h(x),u)) = \nabla V(z)f(z,u) + p(V(z))\nabla V(z)L(h(z)-h(x)) + W(z) \leq$$
$$+ (1-p(V(z)))\nabla V(z)f(z,u) + (1-p(V(z)))W(z) \quad \text{(A.3)}$$
$$+ (1-c)|\nabla V(z)|^2 p(V(z))\frac{(z-x)'Q\big(f(z,u)+L(h(z)-h(x))-f(x,u)\big)}{\nabla V(z)Q(z-x)}$$

Using (A.3), (2.1) and the fact that $0\leq p(V(z))\leq 1$, we obtain:

$$-\frac{\varphi(z,h(x),u))\nabla V(z)Q(z-x)}{|\nabla V(z)|^2} \leq -\frac{1-p(V(z))}{|\nabla V(z)|^2}\nabla V(z)Q(z-x)\big(\nabla V(z)f(z,u)+W(z)\big)$$
$$-(1-c)p(V(z))(z-x)'Q\big(f(z,u)+L(h(z)-h(x))-f(x,u)\big) \leq -(1-c)(z-x)'Q\big(f(z,u)+L(h(z)-h(x))-f(x,u)\big)$$

Combining (2.3), (A.2) and the above inequality, we conclude that (A.1) holds.

Consider a solution of (1.1), (1.2), (2.10), (2.11), (2.12), corresponding to (arbitrary) input $u\in L^\infty([-r-\tau,+\infty);U)$ and satisfying $x(t-r)\in S_1$, $z(t)\in S_2$ for all $t\geq \tau_l$, where $l>0$ is an integer with $\tau_l\geq r$. Next consider the evolution of the mapping $t\to (z(t)-x(t-r))'Q(z(t)-x(t-r))$. Inequality (A.1) and equations (1.1), (2.10) imply that the following inequality holds for $t\geq \tau_l$ a.e.:

$$\frac{d}{dt}\Big((z(t)-x(t-r))'Q(z(t)-x(t-r))\Big) \leq -2c\omega|z(t)-x(t-r)|^2 + 2G_2|Q|\|z(t)-x(t-r)\||w(t)-h(x(t-r))| \quad \text{(A.4)}$$

where $G_2:=\sup\left\{\frac{|\hat{k}(z,y,u)-\hat{k}(z,w,u)|}{|y-w|}: y,w\in\Re^k, z\in S_2, u\in U, y\neq w\right\}$. By virtue of definitions (2.5), (2.6), (2.7), it follows that the constant $G_2$ is well-defined. Since $Q\in\Re^{n\times n}$ is a positive definite matrix there exists a constant $0<K_2\leq |Q|$ with $K_2|x|^2\leq x'Qx$ for all $x\in\Re^n$. Completing the squares and integrating (A.4), we obtain the following estimate for $t\geq \tau_l$:

$$|z(t)-x(t-r)| \leq \exp\left(-\frac{c\omega}{2|Q|}(t-\tau_l)\right)\sqrt{\frac{|Q|}{K_2}}|z(\tau_l)-x(\tau_l-r)| + \sqrt{\frac{2|Q|}{K_2}}\frac{G_2|Q|}{c\omega}\sup_{t_0\leq s\leq t}\left(\exp\left(-\frac{c\omega}{4|Q|}(t-s)\right)|w(s)-h(x(s-r))|\right) \quad \text{(A.5)}$$

Selecting $\sigma>0$ so that $\sigma\leq c\omega/(4|Q|)$, we obtain from (A.5) for $t\geq \tau_l$:

$$\sup_{\tau_l\leq s\leq t}\big(\exp(\sigma s)|z(s)-x(s-r)|\big) \leq \sqrt{\frac{|Q|}{K_2}}|z(t_0)-x(t_0-r)|\exp(\sigma\tau_l) + \sqrt{\frac{2|Q|}{K_2}}\frac{G_2|Q|}{c\omega}\sup_{\tau_l\leq s\leq t}\big(\exp(\sigma s)|w(s)-h(x(s-r))|\big) \quad \text{(A.6)}$$

Finally, notice that since $\sup_{i \geq 0}(\tau_{i+1} - \tau_i) \leq T_s$, it follows from (1.1), (2.11), (2.12) that the following estimate holds for every $t \in [\tau_i, \tau_{i+1})$ with $i \geq l$:

$$|w(t) - h(x(t-r))| \leq T_s G_1 \sup_{\tau_i \leq s \leq t} |z(s) - x(s-r)| \quad (A.7)$$

where $G_1 := \sup\left\{\frac{|\nabla h(x) f(x,u) - \nabla h(z) f(z,u)|}{|x-z|} : x \in S_1, z \in S_2, u \in U, x \neq z\right\}$. Using the inequality $t \leq \tau_i + T_s$ and (A.7) we obtain for all $t \geq \tau_l$:

$$\sup_{\tau_l \leq s \leq t}\left(\exp(\sigma s)|w(s) - h(x(s-r))|\right) \leq T_s G_1 \exp(\sigma T_s) \sup_{\tau_l \leq s \leq t}\left(\exp(\sigma s)|z(s) - x(s-r)|\right) \quad (A.8)$$

Combining (A.6) and (A.8) we get for all $t \geq \tau_l$:

$$\sup_{\tau_l \leq s \leq t}\left(\exp(\sigma s)|z(s) - x(s-r)|\right) \leq \sqrt{\frac{|Q|}{K_2}} \exp(\sigma \tau_l)|z(\tau_l) - x(\tau_l - r)| + T_s G_1 \exp(\sigma T_s) \sqrt{\frac{2|Q|}{K_2}} \frac{G_2|Q|}{c\omega} \sup_{\tau_l \leq s \leq t}\left(\exp(\sigma s)|z(s) - x(s-r)|\right) \quad (A.9)$$

Selecting $T_s > 0$ so that $T_s G_1 \exp(\sigma T_s) \sqrt{\frac{2|Q|}{K_2}} \frac{G_2|Q|}{c\omega} < 1$, we conclude from (A.9) that the following estimate holds:

$$\sup_{\tau_l \leq s \leq t}\left(\exp(\sigma s)|z(s) - x(s-r)|\right) \leq \frac{c\omega\sqrt{|Q|}\exp(\sigma \tau_l)}{c\omega\sqrt{K_2} - T_s G_1 G_2 |Q| \exp(\sigma T_s)\sqrt{2|Q|}} |z(\tau_l) - x(\tau_l - r)|, \quad \forall t \geq \tau_l \quad (A.10)$$

Inequality (3.2) is a direct consequence of (A.10). The proof is complete. ◁

**Proof of Lemma 3.2:** Define:

$$L_X := \sup\left\{\frac{|f(x,u) - f(z,u)|}{|x-z|} : x, z \in S_2, x \neq z, u \in U\right\} \quad (A.11)$$

$$L_U := \sup\left\{\frac{|f(x,u) - f(x,v)|}{|u-v|} : u, v \in U, u \neq v, x \in S_2\right\} \quad (A.12)$$

$$C := \sup\left\{\frac{|\nabla P(x)|}{|x|} : x \in S_1 \setminus \{0\}\right\} \quad (A.13)$$

$$K := \sup\left\{\frac{|k(x) - k(z)|}{|x-z|} : x, z \in S_2, x \neq z\right\} \quad (A.14)$$

Using (2.2) and definitions (3.1), (A.12), (A.13), (A.14), we obtain for all $(x, x_0, u) \in S_1 \times S_1 \times U$:

$$\nabla P(x) f(x,u) \leq -2\mu|x|^2 + |x|CL_U|u - k(x_0)| + |x|CL_U K|x - x_0| \quad (A.15)$$

Consider a solution of (1.1), (1.2), corresponding to (arbitrary) input $u \in L^\infty([-\tau, +\infty); U)$ and satisfying $x(t) \in S_1$ for all $t \geq \tau + jT_H$, where $j \geq 0$ is an integer. Using (1.1) and definitions (3.1), (A.11), (A.12), (A.14), we obtain for all $i \geq j$ and $t \in [iT_H, (i+1)T_H)$:

$$|x(t+\tau) - x(iT_H + \tau)| \leq \int_{iT_H}^{t} |f(x(s+\tau), u(s))| ds \leq \int_{iT_H}^{t} |f(x(s+\tau), k(x(iT_H + \tau)))| ds + T_H L_U \sup_{iT_H \leq s \leq t}\left(|u(s) - k(x(iT_H + \tau))|\right)$$

$$\leq T_H (L_X + L_U K)|x(iT_H + \tau)| + T_H L_X \max_{iT_H \leq s \leq t}\left(|x(s+\tau) - x(iT_H + \tau)|\right) + T_H L_U \sup_{iT_H \leq s \leq t}\left(|u(s) - k(x(iT_H + \tau))|\right) \quad (A.16)$$

Using (A.16), we obtain for all $i \geq j$ for sufficiently small $T_H > 0$ (so that $T_H L_X < 1$):

$$\max_{iT_H \leq s \leq (i+1)T_H}\left(|x(s+\tau) - x(iT_H + \tau)|\right) \leq \frac{(L_X + L_U K)T_H}{1 - T_H L_X}|x(iT_H + \tau)| + \frac{L_U T_H}{1 - T_H L_X} \sup_{iT_H \leq s \leq (i+1)T_H}\left(|u(s) - k(x(iT_H + \tau))|\right) \quad (A.17)$$

Using (A.17), the triangle inequality and a standard causality argument, we obtain for all $i \geq j$ and $t \in [iT_H, (i+1)T_H)$ for sufficiently small $T_H > 0$ (so that $\frac{(L_X + L_U K)T_H}{1 - T_H L_X} < 1$):

$$|x(t+\tau) - x(iT_H + \tau)| \leq \frac{(L_X + L_U K)T_H}{1 - T_H L_X}|x(t+\tau) - x(iT_H + \tau)| + \frac{(L_X + L_U K)T_H}{1 - T_H L_X}|x(t+\tau)| + \frac{L_U T_H}{1 - T_H L_X}\sup_{iT_H \leq s \leq t}(|u(s) - k(x(iT_H + \tau))|)$$

which directly implies the following estimate for all $\sigma > 0$:

$$|x(t+\tau) - x(iT_H + \tau)|\exp(\sigma t) \leq \frac{(L_X + L_U K)T_H \exp(\sigma t)}{1 - (2L_X + L_U K)T_H}|x(t+\tau)|$$
$$+ \frac{L_U T_H \exp(\sigma T_H)}{1 - (2L_X + L_U K)T_H}\sup_{iT_H \leq s \leq t}(|u(s) - k(x(iT_H + \tau))|\exp(\sigma s))$$
(A.18)

Next consider the evolution of the mapping $t \to P(x(t+\tau))$. Inequalities (A.15), (A.18) and equation (1.1) imply that the following differential inequality holds for all $i \geq j$ and $t \in [iT_H, (i+1)T_H)$ a.e.:

$$\frac{d}{dt}P(x(t+\tau)) \leq -\left(2\mu - \frac{CL_U K(L_X + L_U K)T_H}{1 - (2L_X + L_U K)T_H}\right)|x(t+\tau)|^2 + \Theta|x(t+\tau)|\exp(-\sigma t)\sup_{iT_H \leq s \leq t}(|u(s) - k(x(iT_H + \tau))|\exp(\sigma s))$$
(A.19)

where $\Theta := CL_U\left(\frac{KL_U T_H \exp(\sigma T_H)}{1 - (2L_X + L_U K)T_H} + 1\right)$. Completing the squares in (A.19), we get for $t \geq jT_H$ a.e.:

$$\frac{d}{dt}P(x(t+\tau)) \leq -\left(\mu - \frac{CL_U K(L_X + L_U K)T_H}{1 - (2L_X + L_U K)T_H}\right)|x(t+\tau)|^2$$
$$+ \frac{\Theta^2}{4\mu}\exp(-2\sigma t)\sup_{jT_H \leq s \leq t}\left(\left|u(s) - k\left(x\left(\tau + \left[\frac{s}{T_H}\right]T_H\right)\right)\right|^2 \exp(2\sigma s)\right)$$
(A.20)

Since $P \in C^2(\Re^n; \Re_+)$ and since $S_1 \subset \Re^n$ is compact, it follows that there exists $\widetilde{P} > 0$ such that $P(x) \leq \widetilde{P}|x|^2$ for all $x \in S_1$. Selecting $\sigma > 0$, $T_H > 0$ to be sufficiently small (so that $\mu \geq \frac{CL_U K(L_X + L_U K)T_H}{1 - (2L_X + L_U K)T_H} + 4\sigma\widetilde{P}$) and integrating (A.20) we get for all $t \geq jT_H$:

$$P(x(t+\tau)) \leq \exp(-4\sigma(t - jT_H))\widetilde{P}|x(jT_H + \tau)|^2$$
$$+ \frac{\Theta^2}{8\mu\sigma}\exp(-2\sigma t)\sup_{jT_H \leq s \leq t}\left(\left|u(s) - k\left(x\left(\tau + \left[\frac{s}{T_H}\right]T_H\right)\right)\right|^2 \exp(2\sigma s)\right)$$
(A.21)

Using the fact that there exists a constant $K_1 > 0$ with $K_1|x|^2 \leq P(x)$ for all $x \in S_1$ in conjunction with (A.21), we obtain for all $t \geq jT_H$:

$$|x(t+\tau)|\exp(\sigma t) \leq \sqrt{\frac{\widetilde{P}}{K_1}}\exp(\sigma jT_H)x(jT_H + \tau) + \frac{\Theta}{2\sqrt{2\mu\sigma K_1}}\sup_{jT_H \leq s \leq t}\left(\left|u(s) - k\left(x\left(\tau + \left[\frac{s}{T_H}\right]T_H\right)\right)\right|\exp(\sigma s)\right) \quad (A.22)$$

Inequality (3.3) is a direct consequence of estimate (A.22). The proof is complete. ◁

**Proof of Lemma 3.3:** Lemma 2.1 in conjunction with the fact that $x(-r) = x_0 \in S_2$ and definition (3.1) implies that $x(t) \in S_2$ for all $t \in [-r, \tau]$. Applying the inequality $|x(t)| \leq |x_0| + \int_{-r}^{t}|f(x(s), u(s-\tau))|ds$ for the solution $x(t)$ of (1.1) with initial condition $x(-r) = x_0$, corresponding to (arbitrary) input $u \in L^\infty([-r-\tau, 0); U)$ and using (A.11), (A.12) and the Gronwall-Bellman Lemma, we obtain:

$$|x(t)| \leq \exp(L(t+r))(|x_0| + \|u\|), \text{ for all } t \in [-r, \tau] \tag{A.23}$$

where $\|u\| = \sup_{-r-\tau \leq s < 0}(|u(s)|)$ and $L := \max(L_X, L_U)$. Next define:

$$\bar{a} := \max\{|f(x,u)| : (x,u) \in S_2 \times U\} \tag{A.24}$$
$$\Omega := \{\xi \in \Re^n : |\xi - x| \leq \bar{a}(r+\tau), V(x) \leq R\} \tag{A.25}$$

Notice that by virtue of assumption (H1) the set $\Omega$ is compact. We select $\bar{h} > 0$ so that:

$$R + \bar{h}\bar{a}\max\{|\nabla V(x)| : x \in \Omega\} \leq b \tag{A.26}$$

and we select $\bar{P} > 0$ so that:

$$\overline{P} \geq \max\left\{\left|\nabla^2 V(\xi)\right| : |\xi - x| \leq \overline{a}\min(r+\tau,\overline{h}), x \in S_2\right\} \tag{A.27}$$

We next make the following claim.

**Claim 1:** *If $R \leq V(x_i) \leq b$ and $h \leq \min\left(r+\tau,\overline{h}, \dfrac{2}{\overline{a}^2\overline{P}}\min\{W(x): R \leq V(x) \leq b\}\right)$, then*

$$V(x_{i+1}) \leq V(x_i) \tag{A.28}$$

**Proof of Claim 1:** Define the function:
$$g(\lambda) = V(x_i + \lambda(x_{i+1} - x_i)) \tag{A.29}$$
for $\lambda \in [0,1]$. The following equalities hold for all $\lambda \in [0,1]$:

$$\frac{dg}{d\lambda}(\lambda) = \nabla V(x_i + \lambda(x_{i+1} - x_i))(x_{i+1} - x_i) \quad,\quad \frac{d^2g}{d\lambda^2}(\lambda) = (x_{i+1} - x_i)'\nabla^2 V(x_i + \lambda(x_{i+1} - x_i))(x_{i+1} - x_i) \tag{A.30}$$

Moreover, notice that by virtue of (3.1), (A.24) and (2.9), it holds that $|x_{i+1} - x_i| \leq \overline{a}h$. The previous inequality in conjunction with (A.27), (A.30) and the fact that $h \leq \min(r+\tau,\overline{h})$ gives:

$$\left|\frac{d^2g}{d\lambda^2}(\lambda)\right| \leq \overline{a}^2\overline{P}h^2, \text{ for all } \lambda \in [0,1] \tag{A.31}$$

Furthermore, inequality (2.1) in conjunction with (2.9) and (A.30) gives:

$$\frac{dg}{d\lambda}(0) = \nabla V(x_i)\int_{ih}^{(i+1)h} f(x_i, u(s))ds \leq -hW(x_i) \tag{A.32}$$

Combining (A.29), (A.31) and (A.32), we get:
$$V(x_{i+1}) = g(1) \leq V(x_i) - hW(x_i) + \overline{a}^2 h^2 \overline{P}/2 \tag{A.33}$$

Inequality (A.28) is a consequence of (A.33) and the fact that $h \leq \dfrac{2}{\overline{a}^2\overline{P}}\min\{W(x): R \leq V(x) \leq b\}$. ◁

A direct consequence of Claim 1 is that (3.5) holds when $h \leq \min\left(r+\tau,\overline{h}, \dfrac{2}{\overline{a}^2\overline{P}}\min\{W(x): R \leq V(x) \leq b\}\right)$. Indeed, if $R \leq V(x_i) \leq b$ for certain $i = 0,\ldots,N-1$ then the fact that $x_{i+1} \in S_2$ follows from Claim 1 and (3.1). If $V(x_i) \leq R$ for certain $i = 0,\ldots,N-1$, then:

$$V(x_{i+1}) = V(x_i) + \int_0^1 \nabla V(x_i + \lambda(x_{i+1} - x_i))(x_{i+1} - x_i)d\lambda$$

which combined with the fact that $|x_{i+1} - x_i| \leq \overline{a}h$, the fact that $h \leq \min(r+\tau,\overline{h})$ and (A.25) gives:

$$|V(x_{i+1})| \leq R + \overline{a}\overline{h}\max\{|\nabla V(x)| : x \in \Omega\}$$

The above inequality in conjunction with (A.26) and definition (3.1) implies that $x_{i+1} \in S_2$.

We next make the following claim.

**Claim 2:** *Define $e_i := x_i - x(ih - r)$, $i \in \{0,\ldots,N\}$, where $x(t)$ is the solution of (1.1) with initial condition $x(-r) = x_0$ corresponding to input $u \in L^\infty([-r-\tau,0); U)$ and suppose that $h \leq \min\left(r+\tau,\overline{h}, \dfrac{2}{\overline{a}^2\overline{P}}\min\{W(x): R \leq V(x) \leq b\}\right)$. Then*

$$|e_i| \leq \frac{h^2}{2}L_X L(1 + \exp(L(r+\tau)))(|x_0| + \|u\|)\frac{\exp(ihL_X) - 1}{\exp(hL_X) - 1}, \text{ for all } i \in \{1,\ldots,N\} \tag{A.34}$$

*where $L_X \geq 0$ is the Lispchitz constant defined in (A.11).*

**Proof of Claim 2:** Notice that, by virtue of (2.9), the following equation holds for all $i \in \{0,\ldots,N-1\}$:

$$e_{i+1} = e_i + \int_{ih}^{(i+1)h}(f(x_i, u(s-r-\tau)) - f(x(s-r), u(s-r-\tau)))ds \tag{A.35}$$

Notice that Lemma 2.1 in conjunction with the fact that $x_0 \in S_2$ implies that $x(t) \in S_2$ for all $t \in [-r, \tau]$. Hence, definition (A.11) in conjunction with (3.5) implies the following inequality for all $i \in \{0, ..., N-1\}$ and $s \in [ih, (i+1)h]$:

$$|f(x_i, u(s-r-\tau)) - f(x(s-r), u(s-r-\tau))| \leq L_X |x_i - x(s-r)| \quad (A.36)$$

Using the definition $e_i := x_i - x(ih - r)$, definitions (A.11), (A.12) and inequality (A.23), in conjunction with $x(t) \in S_2$ for all $t \in [-r, \tau]$, we get for all $i \in \{0, ..., N-1\}$ and $s \in [ih, (i+1)h]$:

$$|x_i - x(s-r)| \leq |e_i| + |x(s-r) - x(ih-r)| \leq |e_i| + L(s-ih) \max_{s \in [ih,(i+1)h]} (|x(s-r)|) + L(s-ih)\|u\|$$
$$\leq |e_i| + L(s-ih)(1 + \exp(L(r+\tau)))(|x_0| + \|u\|) \quad (A.37)$$

where $L := \max(L_X, L_U)$. Exploiting (A.35), (A.36), (A.37), we obtain for all $i \in \{0, ..., N-1\}$:

$$|e_{i+1}| \leq (1 + hL_X)|e_i| + \frac{h^2}{2} L_X L(1 + \exp(L(r+\tau)))(|x_0| + \|u\|) \leq \exp(hL_X)|e_i| + \frac{h^2}{2} L_X L(1 + \exp(L(r+\tau)))(|x_0| + \|u\|) \quad (A.38)$$

Using the fact $e_0 = 0$, in conjunction with inequality (A.38), gives the desired inequality (A.34). ◁

The desired inequality (3.4) is a direct consequence of estimate (A.34) with $i = N$, the fact that $h := (\tau + r)/N$ and the fact that $\exp(hL_X) - 1 \geq hL_X$ and definitions

$$M_4 := \frac{r+\tau}{2} L(1 + \exp(L(r+\tau)))(\exp((r+\tau)L_X) - 1), \quad N^* = 1 + \left\lceil \frac{r+\tau}{\min\left(r+\tau, \bar{h}, \frac{2}{a^2 \bar{P}} \min\{W(x) : R \leq V(x) \leq b\}\right)} \right\rceil. \quad \triangleleft$$

**Proof of Lemma 3.4:** Let $\sigma > 0$, $T_s > 0$, $T_H > 0$ be sufficiently small constants so that Lemma 3.1 and Lemma 3.2 hold. Let $N \geq N^*$ be an integer, where $N^*$ is the integer constant in Lemma 3.3. Since $z(t) \in S_2$ for all $t \geq jT_H$, it follows from Lemma 3.3 and (3.5) that $\Phi_N(z(q(s)), \breve{u}_{q(s)}) \in S_2$ for all $s \geq jT_H$ and $N \geq N^*$, where $q(s) := [s/T_H]T_H$. Using (A.14) and (2.13) we obtain for all $t \geq jT_H$:

$$\sup_{jT_H \leq s \leq t} (|u(s) - k(x(\tau + q(s)))| \exp(\sigma s)) \leq K \sup_{jT_H \leq s \leq t} (|\Phi_N(z(q(s)), \breve{u}_{q(s)}) - x(\tau + q(s))| \exp(\sigma s)) \quad (A.39)$$

Let $\phi(x_0, u)$ denote the solution of (1.1) at $t = \tau$ with initial condition $x(-r) = x_0$, corresponding to (arbitrary) input $u \in L^\infty([-r-\tau, 0); U)$. It follows that $x(\tau + q(s)) = \phi(x(q(s) - r), \breve{u}_{q(s)})$ for all $s \geq jT_H$. Moreover, using (A.11), Gronwall's Lemma, the fact that all solutions of (1.1) starting from $S_2$ remain in $S_2$ for all times (a consequence of Lemma 2.1) and the fact that $x(t - r - T_H) \in S_1$ for all $t \geq jT_H$, we get:

$$|\phi(z(q(s)), \breve{u}_{q(s)}) - \phi(x(q(s) - r), \breve{u}_{q(s)})| \leq \exp(L_X(r+\tau))|z(q(s)) - x(q(s) - r)| \quad (A.40)$$

Furthermore, Lemma 3.3 implies that there exists $M_4 > 0$ such that for $N \geq N^*$ and $s \geq jT_H$ it holds that:

$$|\phi(z(q(s)), \breve{u}_{q(s)}) - \Phi_N(z(q(s)), \breve{u}_{q(s)})| \leq M_4(|z(q(s))| + \|\breve{u}_{q(s)}\|)/N \quad (A.41)$$

Combining (A.39), (A.40), (A.41) we obtain for all $t \geq jT_H$:

$$\sup_{jT_H \leq s \leq t} (|u(s) - k(x(\tau + q(s)))| \exp(\sigma s)) \leq K\left(\exp(L_X(r+\tau)) + \frac{M_4}{N}\right) \sup_{jT_H \leq s \leq t} (|z(q(s)) - x(q(s) - r)| \exp(\sigma s))$$
$$+ K \frac{M_4}{N} \sup_{jT_H \leq s \leq t} (|x(q(s) - r)| \exp(\sigma s)) + K \frac{M_4}{N} \sup_{jT_H \leq s \leq t} (\|\breve{u}_{q(s)}\| \exp(\sigma s)) \quad (A.42)$$

Moreover, using (A.14), the fact that $k(0) = 0$, the fact that $0 \leq s - q(s) \leq T_H$ and the fact that $x(t - r - T_H) \in S_1$ for all $t \geq jT_H$, we obtain for all $s \geq jT_H$:

$$\|\breve{u}_{q(s)}\| \exp(\sigma s) = \exp(\sigma s) \sup_{-r-\tau \leq \theta < 0} (|u(q(s) + \theta)|)$$
$$\leq \exp(\sigma(s - q(s) + \tau + r)) \sup_{-r-\tau \leq \theta < 0} (|u(q(s) + \theta) - k(x(\tau + q(q(s) + \theta)))| \exp(\sigma(q(s) + \theta))) \quad (A.43)$$
$$+ K \exp(\sigma(s - q(s) - r - \tau)) \sup_{-r-\tau \leq \theta < 0} (|x(\tau + q(q(s) + \theta))| \exp(\sigma q(q(s) + \theta)))$$

Using the fact that $0 \leq s - q(s) \leq T_H$ in conjunction with (A.42), (A.43), we get for all $t \geq jT_H$:

$$\sup_{jT_H \leq s \leq t}\left(|u(s) - k(x(\tau + q(s)))|\exp(\sigma s)\right) \leq K\left(\exp(L_X(r+\tau)) + \frac{M_4}{N}\right)\exp(\sigma T_H)\sup_{jT_H \leq s \leq t}\left(|z(s) - x(s-r)|\exp(\sigma s)\right)$$

$$+ K\frac{M_4}{N}\exp(\sigma(T_H + r + \tau))\sup_{jT_H \leq s \leq t}\left(|x(\tau + q(s) - r - \tau)|\exp(\sigma(q(s) - r - \tau))\right)$$

$$+ K\frac{M_4}{N}\exp(\sigma(2T_H + r + \tau))\left(\sup_{jT_H - r - \tau \leq s \leq t}\left(|u(s) - k(x(\tau + q(s)))|\exp(\sigma s)\right) + K^2 \sup_{(j-1)T_H - r - \tau \leq s \leq t}\left(|x(\tau + s)|\exp(\sigma s)\right)\right)$$

which directly implies for all $t \geq jT_H$:

$$\sup_{jT_H \leq s \leq t}\left(|u(s) - k(x(\tau + q(s)))|\exp(\sigma s)\right) \leq K\left(\exp(L_X(r+\tau)) + \frac{M_4}{N}\right)\exp(\sigma T_H)\sup_{jT_H \leq s \leq t}\left(|z(s) - x(s-r)|\exp(\sigma s)\right)$$

$$+ K\frac{M_4}{N}\exp(\sigma(2T_H + r + \tau))\sup_{jT_H - r - \tau \leq s \leq t}\left(|u(s) - k(x(\tau + q(s)))|\exp(\sigma s)\right) \quad (A.44)$$

$$+ K(1+K)\frac{M_4}{N}\exp(\sigma(2T_H + r + \tau))\sup_{(j-1)T_H - r - \tau \leq s \leq t}\left(|x(\tau + s)|\exp(\sigma s)\right)$$

Selecting $N \geq N^*$ so that $N > KM_4\exp(\sigma(2T_H + r + \tau))$, we get from (A.44) for all $t \geq jT_H$:

$$\sup_{jT_H \leq s \leq t}\left(|u(s) - k(x(\tau + q(s)))|\exp(\sigma s)\right) \leq \frac{K(N\exp(L_X(r+\tau)) + M_4)\exp(2\sigma T_H)}{N - KM_4\exp(\sigma(2T_H + r + \tau))}\sup_{jT_H \leq s \leq t}\left(|z(s) - x(s-r)|\exp(\sigma s)\right)$$

$$+ \frac{KM_4\exp(\sigma(2T_H + r + \tau))}{N - KM_4\exp(\sigma(2T_H + r + \tau))}\sup_{jT_H - r - \tau \leq s \leq jT_H}\left(|u(s) - k(x(\tau + q(s)))|\exp(\sigma s)\right) \quad (A.45)$$

$$+ \frac{K(1+K)M_4\exp(\sigma(2T_H + r + \tau))}{N - KM_4\exp(\sigma(2T_H + r + \tau))}\sup_{(j-1)T_H - r - \tau \leq s \leq t}\left(|x(\tau + s)|\exp(\sigma s)\right)$$

Let $l > 0$ is an integer with $\tau_l \geq r$ and $\tau_l \geq jT_H$. Notice that $\tau_l \leq \max(jT_H, r) + T_s$. It follows from Lemma 3.1 that there exists a constant $M_1 \geq 1$ so that the following inequality holds for all $t \geq jT_H$:

$$\sup_{jT_H \leq s \leq t}\left(\exp(\sigma s)|z(s) - x(s-r)|\right) \leq M_1 \sup_{jT_H \leq s \leq \max(r, jT_H) + T_s}\left(\exp(\sigma s)|z(s) - x(s-r)|\right) \quad (A.46)$$

Combining (A.45) and (A.46) we get for all $t \geq jT_H$:

$$\sup_{jT_H \leq s \leq t}\left(|u(s) - k(x(\tau + q(s)))|\exp(\sigma s)\right) \leq \frac{K(N\exp(L_X(r+\tau)) + M_4)\exp(2\sigma T_H)}{N - KM_4\exp(\sigma(2T_H + r + \tau))}M_1\sup_{jT_H \leq s \leq \max(r, jT_H) + T_s}\left(|z(s) - x(s-r)|\exp(\sigma s)\right)$$

$$+ \frac{KM_4\exp(\sigma(2T_H + r + \tau))}{N - KM_4\exp(\sigma(2T_H + r + \tau))}\left(\sup_{jT_H - r - \tau \leq s \leq jT_H}\left(|u(s) - k(x(\tau + q(s)))|\exp(\sigma s)\right) + (1+K)\sup_{(j-1)T_H - r - \tau \leq s \leq t}\left(|x(\tau + s)|\exp(\sigma s)\right)\right)$$

(A.47)

By virtue of Lemma 3.2, there exist constants $M_2, M_3 > 0$ such that (3.3) holds. Combining (3.3) and (A.47) we get for all $t \geq jT_H$:

$$\sup_{jT_H \leq s \leq t}\left(|x(\tau + s)|\exp(\sigma s)\right) \leq M_2|x(\tau + jT_H)|\exp(\sigma jT_H)$$

$$+ \frac{K(N\exp(L_X(r+\tau)) + M_4)\exp(2\sigma T_H)}{N - KM_4\exp(\sigma(2T_H + r + \tau))}M_1 M_3 \sup_{jT_H \leq s \leq \max(r, jT_H) + T_s}\left(|z(s) - x(s-r)|\exp(\sigma s)\right)$$

$$+ \frac{KM_3 M_4\exp(\sigma(2T_H + r + \tau))}{N - KM_4\exp(\sigma(2T_H + r + \tau))}\sup_{jT_H - r - \tau \leq s \leq jT_H}\left(|u(s) - k(x(\tau + q(s)))|\exp(\sigma s)\right) \quad (A.48)$$

$$+ M_3 \frac{K(1+K)M_4\exp(\sigma(2T_H + r + \tau))}{N - KM_4\exp(\sigma(2T_H + r + \tau))}\sup_{(j-1)T_H - r - \tau \leq s \leq t}\left(|x(\tau + s)|\exp(\sigma s)\right)$$

Selecting $N \geq N^*$ sufficiently large so that $M_3 \frac{K(1+K)M_4\exp(\sigma(2T_H + r + \tau))}{N - KM_4\exp(\sigma(2T_H + r + \tau))} < 1$, we obtain from (A.48) the existence of a constant $\Lambda_1 \geq 1$ for which the following inequality holds for all $t \geq 0$:

$$\sup_{-r-\tau \leq s \leq t}\left(|x(\tau + s)|\exp(\sigma s)\right) \leq \Lambda_1 \sup_{jT_H \leq s \leq \max(r, jT_H) + T_s}\left(|z(s) - x(s-r)|\exp(\sigma s)\right)$$

$$+ \Lambda_1 \sup_{jT_H - r - \tau \leq s \leq jT_H}\left(|u(s) - k(x(\tau + q(s)))|\exp(\sigma s)\right) + \Lambda_1 \sup_{-r-\tau \leq s \leq jT_H}\left(|x(\tau + s)|\exp(\sigma s)\right) \quad (A.49)$$

The definition of the norms $\|x_t\|$ and $\|\tilde{u}_t\|$ give for all $t \geq 0$:

$$\sup_{0 \leq s \leq t}\left(\|x_s\|\exp(\sigma s)\right) \leq \exp(\sigma(r+\tau))\sup_{-r-\tau \leq s \leq t}\left(|x(\tau+s)|\exp(\sigma s)\right) \tag{A.50}$$

$$\sup_{0 \leq s \leq t}\left(\|\tilde{u}_s\|\exp(\sigma s)\right) \leq \exp(\sigma(r+\tau))\sup_{-r-\tau \leq s \leq t}\left(|u(s)|\exp(\sigma s)\right) \tag{A.51}$$

When $t \leq jT_H$ we have from (A.51) that $\sup_{0 \leq s \leq t}\left(\|\tilde{u}_s\|\exp(\sigma s)\right) \leq \exp(\sigma(r+\tau))\sup_{-r-\tau \leq s \leq jT_H}\left(|u(s)|\exp(\sigma s)\right)$. When $t \geq jT_H$, we obtain from (A.47), (A.49), (A.14), the fact that $k(0)=0$, the fact that $0 \leq s - q(s) \leq T_H$ and the fact that $x(t-r-T_H) \in S_1$ for all $t \geq jT_H$, the existence of a constant $\Lambda_2 \geq 1$ for which the following inequality holds for all $t \geq 0$:

$$\sup_{0 \leq s \leq t}\left(\|\tilde{u}_s\|\exp(\sigma s)\right) \leq \exp(\sigma(r+\tau))\sup_{-r-\tau \leq s \leq t}\left(|u(s)|\exp(\sigma s)\right)$$
$$\leq \exp(\sigma(r+\tau))\sup_{-r-\tau \leq s \leq jT_H}\left(|u(s)|\exp(\sigma s)\right) + \exp(\sigma(r+\tau))\sup_{jT_H \leq s \leq t}\left(|u(s)-k(x(\tau+q(s)))|\exp(\sigma s)\right)$$
$$+ K\exp(\sigma(r+\tau+T_H))\sup_{jT_H \leq s \leq t}\left(|x(\tau+s)|\exp(\sigma s)\right)$$
$$\leq \exp(\sigma(r+\tau))\sup_{-r-\tau \leq s \leq jT_H}\left(|u(s)|\exp(\sigma s)\right) + \Lambda_2 \sup_{jT_H \leq s \leq \max(r,jT_H)+T_s}\left(|z(s)-x(s-r)|\exp(\sigma s)\right)$$
$$+ \Lambda_2 \sup_{jT_H-r-\tau \leq s \leq jT_H}\left(|u(s)-k(x(\tau+q(s)))|\exp(\sigma s)\right) + \Lambda_2 \sup_{-r-\tau \leq s \leq jT_H}\left(|x(\tau+s)|\exp(\sigma s)\right)$$

Combining the two cases ($t \leq jT_H$ and $t \geq jT_H$) and using (A.49), (A.50), we obtain the existence of a constant $\Lambda_3 \geq 1$ for which the following inequality holds for all $t \geq 0$:

$$\sup_{0 \leq s \leq t}\left(\|\tilde{u}_s\|\exp(\sigma s)\right) + \sup_{0 \leq s \leq t}\left(\|x_s\|\exp(\sigma s)\right) \leq \Lambda_3 \sup_{-r-\tau \leq s \leq jT_H}\left(|u(s)|\exp(\sigma s)\right) + \Lambda_3 \sup_{-r-\tau \leq s \leq jT_H}\left(|x(\tau+s)|\exp(\sigma s)\right)$$
$$+ \Lambda_3 \sup_{jT_H-r-\tau \leq s \leq jT_H}\left(|u(s)-k(x(\tau+q(s)))|\exp(\sigma s)\right) + \Lambda_3 \sup_{jT_H \leq s \leq \max(r,jT_H)+T_s}\left(|z(s)-x(s-r)|\exp(\sigma s)\right) \tag{A.52}$$

When $t \leq jT_H$ we have that $\sup_{0 \leq s \leq t}\left(|z(s)|\exp(\sigma s)\right) \leq \sup_{0 \leq s \leq jT_H}\left(|z(s)-x(s-r)|\exp(\sigma s)\right) + \exp(\sigma(r+\tau))\sup_{-\tau-r \leq s \leq jT_H}\left(|x(\tau+s)|\exp(\sigma s)\right)$. When $t \geq jT_H$, we obtain from (A.46):

$$\sup_{0 \leq s \leq t}\left(|z(s)|\exp(\sigma s)\right) \leq \sup_{0 \leq s \leq jT_H}\left(|z(s)-x(s-r)|\exp(\sigma s)\right) + \exp(\sigma(r+\tau))\sup_{-\tau-r \leq s \leq t}\left(|x(\tau+s)|\exp(\sigma s)\right) + \sup_{jT_H \leq s \leq t}\left(|z(s)-x(s-r)|\exp(\sigma s)\right)$$
$$\leq (1+M_1)\sup_{0 \leq s \leq \max(r,jT_H)+T_s}\left(|z(s)-x(s-r)|\exp(\sigma s)\right) + \exp(\sigma(r+\tau))\sup_{-\tau-r \leq s \leq t}\left(|x(\tau+s)|\exp(\sigma s)\right)$$

Combining the two cases ($t \leq jT_H$ and $t \geq jT_H$) and using (A.49), (A.52), we obtain the existence of a constant $\Lambda_4 \geq 1$ for which the following inequality holds for all $t \geq 0$:

$$\sup_{0 \leq s \leq t}\left(\|\tilde{u}_s\|\exp(\sigma s)\right) + \sup_{0 \leq s \leq t}\left(\|x_s\|\exp(\sigma s)\right) + \sup_{0 \leq s \leq t}\left(|z(s)|\exp(\sigma s)\right)$$
$$\leq \Lambda_4 \sup_{-r-\tau \leq s \leq jT_H}\left(|u(s)|\exp(\sigma s)\right) + \Lambda_4 \sup_{0 \leq s \leq \max(r,jT_H)+T_s}\left(|z(s)-x(s-r)|\exp(\sigma s)\right) \tag{A.53}$$
$$+ \Lambda_4 \sup_{jT_H-r-\tau \leq s \leq jT_H}\left(|u(s)-k(x(\tau+q(s)))|\exp(\sigma s)\right) + \Lambda_4 \sup_{-r-\tau \leq s \leq jT_H}\left(|x(\tau+s)|\exp(\sigma s)\right)$$

Moreover, using (A.14), the fact that $k(0)=0$, the fact that $0 \leq s - q(s) \leq T_H$ and the fact that $x(t-r-T_H) \in S_1$ for all $t \geq jT_H$, we obtain from (A.49) for all $t \geq 0$:

$$\sup_{jT_H-r-\tau \leq s \leq jT_H}\left(|u(s)-k(x(\tau+q(s)))|\exp(\sigma s)\right) \leq \sup_{jT_H-r-\tau \leq s \leq jT_H}\left(|u(s)|\exp(\sigma s)\right) + \sup_{jT_H-r-\tau \leq s \leq jT_H}\left(|k(x(\tau+q(s)))|\exp(\sigma s)\right)$$
$$\leq \sup_{-r-\tau \leq s \leq jT_H}\left(|u(s)|\exp(\sigma s)\right) + K\sup_{jT_H-r-\tau \leq s \leq jT_H}\left(|x(\tau+q(s))|\exp(\sigma s)\right)$$
$$\leq \sup_{-r-\tau \leq s \leq jT_H}\left(|u(s)|\exp(\sigma s)\right) + K\exp(\sigma T_H)\sup_{-r-\tau \leq s \leq jT_H}\left(|x(\tau+s)|\exp(\sigma s)\right)$$

The above inequality in conjunction with (A.53) implies the desired inequality (3.6). The proof is complete. ◁